\documentclass[11pt]{article}
\usepackage{amsfonts}

\usepackage[dvips]{graphicx}
\usepackage{amsmath}

\newtheorem{theorem}{Theorem}[subsection]

\newtheorem{corollary}[theorem]{Corollary}

\newtheorem{definition}[theorem]{Definition}
\newtheorem{example}[theorem]{Example}

\newtheorem{lemma}[theorem]{Lemma}

\newtheorem{proposition}[theorem]{Proposition}

\newenvironment{proof}[1][Proof]{\textbf{#1.} }{\ \rule{0.5em}{0.5em}}

\input{tcilatex}

\begin{document}

\title{Crystal bases and combinatorics of infinite rank quantum groups}
\author{C\'{e}dric Lecouvey \\
Laboratoire de Math\'{e}matiques Pures et Appliqu\'{e}es Joseph Liouville\\
B.P. 699 62228 Calais Cedex}
\date{}
\maketitle

\begin{abstract}
The tensor powers of the vector representation associated to an infinite
rank quantum group decompose into irreducible components with multiplicities
independant of the infinite root system considered. Although the irreducible
modules obtained in this way are not of highest weight, they admit a crystal
basis and a canonical basis. This permits in particular to obtain for each
familly of classical Lie algebras a Robinson-Schensted correspondence on
biwords defined on infinite alphabets. We then depict a structure of
bi-crystal on these bi-words. This RS-correspondence yields also a plactic
algebra and plactic Schur functions distinct for each infinite root system.
Surprisingly, the algebras spanned by these plactic Schur functions are all
isomorphic to the algebra of symmetric functions.
\end{abstract}

\section{ Introduction}

The diagonal expansion of the Cauchy kernel as a sum over all partitions $%
\lambda$%
\begin{equation}
\prod_{i\in I}\prod_{j\in J}(1-x_{i}y_{j})^{-1}=\sum_{\lambda}s_{\lambda
}(X)s_{\lambda}(Y)  \label{CI}
\end{equation}
permits to characterize the Schur functions $s_{\lambda}$ (see \cite{LA2}).
This Cauchy identity can be interpreted in terms of the combinatorics of
Young tableaux by using the Robinson-Schensted correspondence for biwords
defined on the totally ordered alphabets $\mathcal{X}=\{x_{i}\mid i\in I\}$
and $\mathcal{Y}=\{y_{j}\mid j\in J\}$ (see \cite{Fu}). Recall that a biword 
$W$ on $\mathcal{X}$ and $\mathcal{Y}$ can be written 
\begin{equation*}
W=\left( 
\begin{array}{lllllllllllllll}
y_{1} & \cdot & \cdot & y_{1} & y_{2} & \cdot & \cdot & y_{2} & \cdot & \cdot
& \cdot & y_{k} & \cdot & \cdot & y_{k} \\ 
b_{1} & \cdot & \cdot & e_{1} & b_{2} & \cdot & \cdot & e_{2} & \cdot & \cdot
& \cdot & b_{k} & \cdot & \cdot & e_{k}
\end{array}
\right) =\left( 
\begin{array}{l}
w_{y} \\ 
w_{x}
\end{array}
\right)
\end{equation*}
where for any $p\in\{1,...,k\}$, $b_{p}\cdot\cdot\cdot e_{p}$ is a word on $%
\mathcal{X}$ such that $b_{p}\geq\cdot\cdot\cdot\geq e_{p}$. The
RS-correspondence associates to each biword $W\in\mathcal{W}$ a pair of
semistandard tableaux $(P,Q)$ with the same shape where $P$ is obtained by
applying the Schensted bumping algorithm to $w_{x}\in\mathcal{X}$ and $Q$ is
the corresponding recording tableau with evaluation $w_{y}.$ By identifying
each word $w$ of length $\ell$ on $\mathcal{X}$ with the biword $\binom
{12\cdot\cdot\cdot l}{w}$, one also obtains a correspondence between words $%
w $ on $\mathcal{X}$ and pairs $(P(w),Q(w))$ where $Q(w)$ is a standard
tableau.\ This makes naturally appear the plactic monoid defined as the
quotient of the free monoid on $\mathcal{X}$ by the Knuth relations 
\begin{equation*}
\left\{ 
\begin{tabular}{l}
$abx=bax$ if $a<x\leq b$ \\ 
$abx=axb$ if $x\leq a<b$%
\end{tabular}
\right. .
\end{equation*}
More precisely, two words $w_{1}$ and $w_{2}$ of length $\ell$ on $\mathcal{X%
}$ are congruent in the plactic monoid if and only if $P(w_{1})=P(w_{2})$.
Thank to this monoid, it is possible to embed the algebra of symmetric
functions into a noncommutative algebra called the plactic algebra which
permits to obtain an illuminating proof of the Littlewood-Richardson rule 
\cite{loth}. The plactic monoid also naturally occur in the theory of
crystal basis developed by Kashiwara \cite{Ka2}. To each finite dimensional $%
U_{q}(\frak{sl}_{n})$-module $M$ is associated a crystal graph whose
decomposition in connected components reflects the decomposition of $M$ into
its irreducible components. Let $B^{\otimes\ell}$ be the crystal graph of
the $\ell$-th tensor power of the vector representation (i.e. the
fundamental $n$-dimensional $U_{q}(\frak{sl}_{n})$-module). The vertices of $%
B^{\otimes\ell}$ can be identified with the words of length $\ell$ on $%
\mathcal{X}$. The plactic monoid is then obtained by identifying the
vertices of $B^{\otimes\ell}$ which occur at the same place in two
isomorphic connected components. Similarly, the Robinson-Schensted
correspondence on biwords can be interpreted as an isomorphism of crystal
graphs (see \ref{subsec_RS}). Moreover it naturally endows $\mathcal{W}$
with the structure of a bi-crystal graph \cite{DC}, \cite{LA2}.

\noindent Crystal basis theory permits also to define plactic monoids and to
obtain RS-type correspondences on \textit{words} for the quantum groups $%
U_{q}(\frak{so}_{2n+1}),U_{q}(\frak{sp}_{2n})$ and $U_{q}(\frak{so}_{2n})$ 
\cite{Lec}, \cite{lec2}. These correspondences are based on insertion
algorithms for tableaux of types $B_{n},C_{n}$ and $D_{n}$ \cite{KN} which
extend the notion of semistandard tableaux to the other classical types. For
each type $\Delta_{n},$ $\Delta=B,C,D,$ one can also consider biwords of
type $\Delta_{n}$ and associate to each such biword $W$ a pair $%
(P_{\Delta_{n}},Q_{\Delta_{n}})$ where $P_{\Delta_{n}}$ is tableau of type $%
\Delta_{n}$ and $Q_{\Delta_{n}}$ a recording tableau. This recording tableau
consists in a sequence of Young diagrams $Q_{1},...,Q_{\ell}$ where the
number of boxes in $Q_{i}$ can be strictly greater than that of $Q_{i+1}.\;$%
This implies in particular that the association $W\mapsto(P_{\Delta_{n}},Q_{%
\Delta_{n}})$ cannot yield a one-to-one correspondence on \textit{biwords}
of type $\Delta_{n}$. When the previous association is restricted to \textit{%
words} $w$ via the identification $w=\binom{12\cdot\cdot\cdot l}{w}$, the
sequence of Young diagrams $Q_{1},...,Q_{\ell}$ obtained is such that $Q_{i}$
and $Q_{i+1}$ differ by at most one box. In this special case, it is
established in \cite{Lec} and \cite{lec2} that the map $w\rightarrow(P_{%
\Delta_{n}},Q_{\Delta_{n}})$ is a one-to-one correspondence.

\bigskip

\noindent Let $\frak{g}_{\infty}$ be a Lie algebra of type $\Delta_{\infty
}=A_{\infty},B_{\infty},C_{\infty}$ or $D_{\infty}$.\ In this paper, we show
that the combinatorial constructions and properties related to the root
systems of type $A_{n},$ $n\geq1$ we have described above can be generalized
to the infinite root systems $\Delta_{\infty}$. We show in particular that
for each infinite root system $\Delta_{\infty}$, there exists a RS-type
correspondence (thus bijective) $W\mapsto(P_{\infty},Q_{\infty})$ on biwords
where $P_{\infty}$ and $Q_{\infty}$ are tableaux respectively of type $%
\Delta_{\infty}$ and semi-standard with the same shape. This yields the
structure of a bi-crystal (that is for $U_{q}(\frak{g}_{\infty})$ and $U_{q}(%
\frak{sl}_{\infty}))$ on biwords. The combinatorial description of these
RS-type correspondences is based on the existence of a plactic monoid for
each type $\Delta_{\infty}$. The defining relations of these monoids depend
on the root system considered. They permit to introduce a plactic algebra
and plactic Schur functions of type $\Delta_{\infty}$. We show that the
plactic Schur functions of type $\Delta_{\infty}$ span a commutative algebra
isomorphic to the algebra of symmetric functions. In particular the four
algebras obtained in this way are isomorphic. Our correspondence on biwords
gives also Cauchy-type identities analogous to (\ref{CI}) for types $%
B_{\infty},C_{\infty}$ and $D_{\infty}$. These formulas hold only in
infinite rank case and do not permit to recover Littelwood's formulas
described in \cite{lit}. To establish these combinatorial results, we need
an extension of the crystal basis theory to certain irreducible $U_{q}(\frak{%
g}_{\infty})$-modules $M_{\infty}$ which are not of highest weight. These
modules are the irreducible components which appear in the decomposition of
the tensor powers of the vector representations associated to the quantum
groups of infinite rank. They are irreducible $U_{q}(\frak{g}_{\infty})$%
-modules $V_{\infty }(\lambda)$ naturally labelled by partitions.\ Since we
have not found references on these modules in the literature, we have
devoted two sections of the paper to the exposition of their algebraic
properties. We establish notably the decomposition 
\begin{equation}
V_{\infty}(\lambda)\otimes V_{\infty}(\lambda)\simeq\bigoplus_{\nu}V_{\nu
}(\lambda)^{\oplus c_{\lambda,\mu}^{\nu}}  \label{deci}
\end{equation}
where $c_{\lambda,\mu}^{\nu}$ is the Littelwood-Richardson coefficient
associated to $\lambda,\mu$ and $\nu$.\ Note that the decomposition (\ref
{deci}) is independent of the infinite root system considered. We also show
how to define crystal bases for the modules $M_{\infty}$ so that the
irreducible modules $V_{\infty}(\lambda)$ have a unique crystal basis (up to
an overall renormalization) and a canonical basis. The general philosophy of
the paper is to introduce the algebraic and combinatorial objects related to
the infinite rank quantum groups as relevant direct limits of their
analogues in finite rank.

\bigskip

\noindent The paper is organized as follows. In section $2,$ we review the
necessary background on quantum infinite rank Lie algebras and recall the
crystal basis theory for their highest weight modules which is essentially
the same as in finite rank. Section $3$ is devoted to the definition of the
modules $V_{\infty}(\lambda)$ and the proof of the Littelwood-Richardson
rule (\ref{deci}).\ We introduce the crystal bases of the module $M_{\infty}$
and the canonical basis of the module $V_{\infty}(\lambda)$ in Section $4$.
The RS-type correspondence and the bi-crystal structures are described in
Section $5$. Finally we prove in Section $6$ that the algebras spanned by
the plactic Schur functions are all isomorphic to the algebra of symmetric
functions and establish Cauchy types identities.

\section{Background\label{secBa}}

In this section we introduce the infinite rank Lie algebras we consider in
this paper together with the corresponding quantum groups.\ We also include
a few considerations on their highest weight modules to avoid a possible
confusion with the modules $V_{\infty}(\lambda)$ (which are not of highest
weight) of Section \ref{sec_modu}. The reader is referred to \cite{Jan}, 
\cite{Kac}, \cite{Kan}, \cite{Ka1} and \cite{Ka2} for more details.

\subsection{Quantum groups of infinite rank}

Let $\frak{gl}_{\pm\infty}$ be the Lie algebra of complex $\mathbb{Z}\times%
\mathbb{Z}$ matrices $A=(a_{i,j})$ with a finite number of nonzero entries.\
The Lie algebra $\frak{gl}_{\pm\infty}$ contains Lie subalgebras $\frak{g}%
_{A_{\infty}},\frak{g}_{B_{\infty}},\frak{g}_{C_{\infty}}$ and $\frak{g}%
_{D_{\infty}}$ of type $A_{\infty},B_{\infty},C_{\infty}$ and $D_{\infty}$
(see \cite{KAC2}, \S\ 7.11). To make our notation homogeneous we write $%
A_{\infty}$ the infinite Dynkin diagram which is denoted by $A_{+\infty}$ in 
\cite{KAC2}$.$

\begin{align}
A_{\infty }& :\overset{1}{\circ }-\overset{2}{\circ }-\overset{3}{\circ }-%
\overset{4}{\circ }-\cdot \cdot \cdot \hspace{0.5cm}B_{\infty }:\overset{0}{%
\circ }\Longleftarrow \overset{1}{\circ }-\overset{2}{\circ }-\overset{3}{%
\circ }-\overset{4}{\circ }-\cdot \cdot \cdot  \label{DD} \\
C_{\infty }& :\overset{0}{\circ }\Longrightarrow \overset{1}{\circ }-%
\overset{2}{\circ }-\overset{3}{\circ }-\overset{4}{\circ }-\cdot \cdot
\cdot \hspace{0.5cm}D_{\infty }: 
\begin{tabular}{l}
$\overset{1}{\circ }$ \\ 
$\ \ \backslash $ \\ 
$\ \ \ \ \overset{2}{\circ }$ \\ 
$\ \ /$ \\ 
$\underset{0}{\circ }$%
\end{tabular}
-\overset{3}{\circ }-\overset{4}{\circ }-\cdot \cdot \cdot  \notag
\end{align}
Note that the infinite rank Lie algebras we use in this paper can all be
realized as Lie algebras of infinite matrices with a finite number of
nonzero entries. The nodes of the Dynkin diagrams (\ref{DD}) are labelled by 
$I$, where $I=\mathbb{N}-\{0\}$ for type $A_{\infty }$ and $I=\mathbb{N}$\
for types $B_{\infty },C_{\infty }$ and $D_{\infty }.$

\noindent We denote by $e_{i},f_{i},t_{i}$ $i\in I,$ the set of Chevalley
generators of $U_{q}(\frak{g}_{\infty }).$ Accordingly to our convention for
the labelling of the infinite Dynkin diagrams (\ref{DD}), the fundamental
weights of $U_{q}(\frak{g}_{\infty })$ belong to $\left( \frac{\mathbb{Z}}{2}%
\right) ^{\infty }$.\ More precisely we have $\mathbf{\omega }%
_{i}=(0^{i},1,1,.....)$ for $i\geq 2$ and also $i=1$ for $\frak{g}_{\infty
}\neq \frak{g}_{D_{\infty }}$, $\mathbf{\omega }_{0}^{C_{\infty
}}=(1,1,1,.....),$ $\mathbf{\omega }_{0}^{B_{\infty }}=\mathbf{\omega }%
_{0}^{D_{\infty }}=(\frac{1}{2},\frac{1}{2},...)$ and $\mathbf{\omega }%
_{1}^{D_{\infty }}=(-\frac{1}{2},\frac{1}{2},...).$ The weight lattice $P$
of $U_{q}(\frak{g}_{\infty })$ can be considered as the $\mathbb{Z}$%
-sublattice of $\left( \frac{\mathbb{Z}}{2}\right) ^{\infty }$ generated by
the $\mathbf{\omega }_{i},$ $i\in I.$ So a weight $\mathbf{\delta }$
associated to $U_{q}(\frak{g}_{\infty })$ is an infinite sequence $\mathbf{%
\delta }=(\delta _{1},\delta _{2},....)$ of integers (resp.\ half integers)
for $\frak{g}_{\infty }=\frak{g}_{A_{\infty }},\frak{g}_{C_{\infty }}$
(resp.\ $\frak{g}_{\infty }=\frak{g}_{B_{\infty }},\frak{g}_{D_{\infty }})$.
We denote by $P^{+}$ the cone of dominant weights of $U_{q}(\frak{g}_{\infty
}),$ that is the set of finite linear combinations of fundamental weights
with nonnegative integer coefficients. They correspond to the infinite
sequences $\mathbf{\lambda }=(\lambda _{1},\lambda _{2},...)\in P$ such that 
\begin{equation*}
0\leq \lambda _{1}\leq \lambda _{2}\leq \lambda _{3}\leq \cdot \cdot \cdot
\leq \lambda _{i}=\lambda _{i+1}=\lambda _{i+2}=\cdot \cdot \cdot 
\end{equation*}
except for type $D_{\infty }$ in which case the inequality $\lambda _{1}\leq
\lambda _{2}$ above is replaced by $\left| \lambda _{1}\right| \leq \lambda
_{2}.$

\noindent Let $\varepsilon _{i}$ be the $i$-th standard basis vector\ of $%
\mathbb{Z}^{\infty }$. The simple roots $\{\alpha _{i}\mid i\in I\}$ are
given by $\alpha _{i}=\varepsilon _{i+1}-\varepsilon _{i}$ for $i\neq 0,$ $%
\alpha _{0}^{B_{\infty }}=-\varepsilon _{1},$ $\alpha _{0}^{C_{\infty
}}=-2\varepsilon _{1}$ and $\alpha _{0}^{D_{\infty }}=-\varepsilon
_{1}-\varepsilon _{2}.$ Write $Q$ for the root lattice and set $Q^{+}=Q\cap
P^{+}.$

\noindent For any nonnegative integer $n$ set $I_{n}=\{i\in I\mid i<n\}$ and
denote by $U_{q}(\frak{g}_{n})$ the subalgebra of $U_{q}(\frak{g}_{\infty})$
generated by $e_{i},f_{i},t_{i},$ $i\in I_{n}.$ Clearly $U_{q}(\frak{g}_{n})$
is the quantum group associated to the finite root system obtained by
considering only the simple roots labelled by $i\in I_{n}$ in (\ref{DD}).
For any $\mathbf{\delta}=(\delta_{1},\delta_{2},....)\in P,$ set $\pi _{n}(%
\mathbf{\delta}\mathbf{)}=(\delta_{1},\delta_{2},...,\delta_{n}).$ Denote by 
$P_{n}$ and $Q_{n}$ the weight and root lattices of $U_{q}(\frak{g}_{n})$.\
One has $P_{n}=\pi_{n}(P)$ and $Q_{n}=\pi_{n}(Q).$ Moreover the set of
dominant weights for $U_{q}(\frak{g}_{n})$ verifies $P_{n}^{+}=%
\pi_{n}(P_{+}) $. For any integer $p\geq n,$ it will be convenient to
identify $P_{n}$ with the sublattice of $P_{p}$ containing the weights of
the form $(\delta _{1},...,\delta_{n},0,...,0).$ Note that the dominant
weights of $P_{n}^{+}$ are not dominant in $P_{p}$ via this identification.
For any $\delta =(\delta_{1},...,\delta_{n})\in P_{n},$ we set $\left|
\delta\right| =\delta_{1}+\cdot\cdot\cdot+\delta_{n}.$

\noindent Given two $U_{q}(\frak{g}_{\infty})$-modules $M$ and $N$, the
structure of $U_{q}(\frak{g}_{\infty})$-module on $M\otimes N$ is defined by
putting 
\begin{gather}
q^{h}(u\otimes v)=q^{h}u\otimes q^{h}v,  \label{tensor1} \\
e_{i}(u\otimes v)=e_{i}u\otimes t_{i}^{-1}v+u\otimes e_{i}v,  \label{tensor2}
\\
f_{i}(u\otimes v)=f_{i}u\otimes v+t_{i}u\otimes f_{i}v.  \label{tensor3}
\end{gather}
for any $i\in I$.

\subsection{Modules of the category $\mathcal{O}_{\mathrm{int}}$}

The reader is referred to \cite{Kan} for a complete exposition of the
results on the category $\mathcal{O}_{\mathrm{int}}$ in finite rank case.\
We also denote by $\mathcal{O}_{\mathrm{int}}$ the category of $U_{q}(\frak{g%
}_{\infty})$-modules $\mathcal{M}$ satisfying the following conditions:

\begin{enumerate}
\item  $\mathcal{M}$ has a weight space decomposition, that is $\mathcal{M}%
=\oplus_{\mathbf{\delta}\in P}\mathcal{M}_{\mathbf{\delta}}$ where $\mathcal{%
M}_{\mathbf{\delta}}=\{v\in\mathcal{M}\mid q^{h_{i}}v=q^{\mathbf{\delta}%
(h_{i})}$ for all $i\in\mathbb{N}\}$ and $\dim \mathcal{M}_{\mathbf{\delta}%
}<\infty.$

\item  There exists a finite number of weights $\mathbf{\lambda}^{(1)},%
\mathbf{\lambda}^{(2)},...,\mathbf{\lambda}^{(s)}$ such that 
\begin{equation*}
\mathcal{M}_{\mathbf{\delta}}\neq\{0\}\Longrightarrow\mathbf{\delta}\in D(%
\mathbf{\lambda}^{(1)})\cup\cdot\cdot\cdot\cup D(\mathbf{\lambda}^{(s)})
\end{equation*}
where for any $k=1,...,s,$ $D(\mathbf{\lambda}^{(k)})=\{\mathbf{\gamma}\in
P\mid\mathbf{\gamma}\leq\mathbf{\lambda}^{(k)}\}.$

\item  The generators $e_{i},f_{i},$ $i\in\mathbb{N}$ are locally nilpotent
on $\mathcal{M}.$
\end{enumerate}

\noindent By using similar arguments to those exposed in \cite{Kan}, we
obtain that the $U_{q}(\frak{g}_{\infty})$-modules of the category $\mathcal{%
O}_{\mathrm{int}}$ are totally reducible and admit crystal bases. If $(L,B)$
and $(L^{\prime},B^{\prime})$ are crystal bases of the $U_{q}(\frak{g}%
_{\infty})$-modules $\mathcal{M}$ and $\mathcal{M}^{\prime}$ belonging to $%
\mathcal{O}_{\mathrm{int}}$, then $(L\otimes L^{\prime},$ $B\otimes
B^{\prime})$ with $B\otimes B^{\prime}=\{b\otimes b^{\prime};$ $b\in
B,b^{\prime}\in B^{\prime}\}$ is a crystal basis of $\mathcal{M}\otimes%
\mathcal{M}^{\prime}$. The action of the Kashiwara operators $\widetilde{e}%
_{i}$ and $\widetilde{f}_{i}$ on $B\otimes B^{\prime}$ is given by:

\begin{align}
\widetilde{f_{i}}(u\otimes v) & =\left\{ 
\begin{tabular}{c}
$\widetilde{f}_{i}(u)\otimes v$ if $\varphi_{i}(u)>\varepsilon_{i}(v)$ \\ 
$u\otimes\widetilde{f}_{i}(v)$ if $\varphi_{i}(u)\leq\varepsilon_{i}(v)$%
\end{tabular}
\right.  \label{TENS1} \\
& \text{and}  \notag \\
\widetilde{e_{i}}(u\otimes v) & =\left\{ 
\begin{tabular}{c}
$u\otimes\widetilde{e_{i}}(v)$ if $\varphi_{i}(u)<\varepsilon_{i}(v)$ \\ 
$\widetilde{e_{i}}(u)\otimes v$ if$\varphi_{i}(u)\geq\varepsilon_{i}(v)$%
\end{tabular}
\right.  \label{TENS2}
\end{align}
where $\varepsilon_{i}(u)=\max\{k;\widetilde{e}_{i}^{k}(u)\neq0\}$ and $%
\varphi_{i}(u)=\max\{k;\widetilde{f}_{i}^{k}(u)\neq0\}$.

\noindent The set $B$ may be endowed with a combinatorial structure called
the crystal graph of $\mathcal{M}$. Crystal graphs for $U_{q}(\frak{g}%
_{\infty})$-modules in $\mathcal{O}_{\mathrm{int}}$ are infinite oriented
colored graphs with colors $i\in I.$ An arrow $a\overset{i}{\rightarrow}b$
means that $\widetilde{f}_{i}(a)=b$ and $\widetilde{e}_{i}(b)=a$. The
decomposition of $\mathcal{M}$ into its irreducible components is reflected
into the decomposition of $B$ into its connected components.

\noindent The weight of $b\in B$ is denoted $\mathrm{wt}(b)$ and we have 
\begin{equation}
\mathrm{wt}(b)=\sum_{i\in I}(\varphi_{i}(b)-\varepsilon_{i}(b))\mathbf{%
\omega }_{i}\in P.  \label{wt}
\end{equation}

\noindent The following lemma is a straightforward consequence of (\ref
{TENS1}) and (\ref{TENS2}).

\begin{lemma}
\label{lem_plu_hp}Let $u\otimes v$ $\in$ $B\otimes B^{\prime}$ $u\otimes v$
is a highest weight vertex of $B\otimes B^{\prime}$ if and only if for any $%
i\in I,$ $\widetilde{e}_{i}(u)=0$ (i.e. $u$ is of highest weight) and $%
\varepsilon_{i}(v)\leq\varphi_{i}(u).$
\end{lemma}

\subsection{The highest weight modules $\mathcal{V}(\mathbf{\protect\lambda }%
)$}

\noindent For any dominant weight $\mathbf{\lambda }\in P_{+}$, we denote by 
$\mathcal{V}(\mathbf{\lambda })$ the irreducible highest weight $U_{q}(\frak{%
g}_{\infty })$-module of highest weight $\mathbf{\lambda .\;}$It is an
infinite dimensional module which belongs to the category $\mathcal{O}_{%
\mathrm{int}}$.\ Thus $\mathcal{V}(\mathbf{\lambda })$ admits a crystal
basis $(L\mathbf{(\lambda ),}B\mathbf{(\lambda )})$ unique up to an overall
scalar factor. More precisely there exists a highest weight vector $v_{%
\mathbf{\lambda }}$ in $\mathcal{V}(\mathbf{\lambda })$ such that 
\begin{align}
L\mathbf{(\lambda )}& =\bigoplus_{r\geq 0,i_{k}\in I}\mathbf{A}(q)\widetilde{%
f}_{i_{1}}\cdot \cdot \cdot \widetilde{f}_{i_{r}}v_{\mathbf{\lambda }}
\label{bc} \\
B\mathbf{(\lambda )}& =\{\widetilde{f}_{i_{1}}\cdot \cdot \cdot \widetilde{f}%
_{i_{r}}v_{\mathbf{\lambda }}+qL\mathbf{(\lambda )}\in L\mathbf{(\lambda )}%
/qL\mathbf{(\lambda )}\mid r\geq 0,i_{k}\in I\}-\{0\}  \notag
\end{align}
where $\mathbf{A}(q)$ is the subalgebra of $\mathbb{Q}(q)$ consisting of the
rational functions without pole at $q=0$ and the $\widetilde{f}_{i_{k}}$'s
are the Kashiwara operators.

\begin{lemma}
\label{lm_tec}Consider $v_{\mathbf{\lambda }}$ and $v_{\mathbf{\lambda }%
}^{\prime }$ two highest weight vectors in $\mathcal{V}($\textbf{$\lambda $}$%
)$ and denote respectively by $(L\mathbf{(}$\textbf{$\lambda $}$\mathbf{),}B%
\mathbf{(}$\textbf{$\lambda $}$\mathbf{)})$ and $(L^{\prime }\mathbf{(}$%
\textbf{$\lambda $}$\mathbf{),}B^{\prime }\mathbf{(}$\textbf{$\lambda $}$%
\mathbf{)})$ the crystal bases obtained from $v_{\mathbf{\lambda }}$ and $v_{%
\mathbf{\lambda }}^{\prime }$ as in (\ref{bc}).\ Then $(L\mathbf{(}$\textbf{$%
\lambda $}$\mathbf{),}B\mathbf{(}$\textbf{$\lambda $}$\mathbf{)})=(L^{\prime
}\mathbf{(}$\textbf{$\lambda $}$\mathbf{),}B^{\prime }\mathbf{(}$\textbf{$%
\lambda $}$\mathbf{)})$ if and only if there exists $K\in \mathbf{A}(q)$
with $K(0)=1$ such that $v_{\mathbf{\lambda }}^{\prime }=K(q)v_{\mathbf{%
\lambda }}.$
\end{lemma}

\begin{proof}
Suppose that $(L\mathbf{(}$\textbf{$\lambda $}$\mathbf{),}B\mathbf{(}$%
\textbf{$\lambda $}$\mathbf{)})=(L^{\prime }\mathbf{(}$\textbf{$\lambda $}$%
\mathbf{),}B^{\prime }\mathbf{(}$\textbf{$\lambda $}$\mathbf{)})$.\ Since $L%
\mathbf{(}$\textbf{$\lambda $}$\mathbf{)}=L^{\prime }\mathbf{(}$\textbf{$%
\lambda $}$\mathbf{)}$ we must have by (\ref{bc}) $v_{\mathbf{\lambda }%
}^{\prime }\in L\mathbf{(}$\textbf{$\lambda $}$\mathbf{)}$.\ Thus there
exists $K(q)\in \mathbf{A}(q)$ such that $v_{\mathbf{\lambda }}^{\prime
}=K(q)v_{\mathbf{\lambda }}.$ Moreover $v_{\mathbf{\lambda }}^{\prime }+qL%
\mathbf{(}$\textbf{$\lambda $}$\mathbf{)}\in B\mathbf{(}$\textbf{$\lambda $}$%
\mathbf{)}$ for $B^{\prime }\mathbf{(}$\textbf{$\lambda $}$\mathbf{)}=B%
\mathbf{(}$\textbf{$\lambda $}$\mathbf{)}$. Since $v_{\mathbf{\lambda }}+qL%
\mathbf{(}$\textbf{$\lambda $}$\mathbf{)}$ is the unique vertex in $B\mathbf{%
(}$\textbf{$\lambda $}$\mathbf{)}$ of weight \textbf{$\lambda $}, we obtain $%
v_{\mathbf{\lambda }}^{\prime }+qL\mathbf{(}$\textbf{$\lambda $}$\mathbf{)=}%
v_{\mathbf{\lambda }}+qL\mathbf{(}$\textbf{$\lambda $}$\mathbf{)}$.\ By
using the equality $v_{\mathbf{\lambda }}^{\prime }=K(q)v_{\mathbf{\lambda }}
$, we derive $K(0)=1$.

\noindent Conversely, if $v_{\mathbf{\lambda }}^{\prime }=K(q)v_{\mathbf{%
\lambda }}$ with $K(0)=1$, we have $L\mathbf{(}$\textbf{$\lambda $}$\mathbf{)%
}=L^{\prime }\mathbf{(}$\textbf{$\lambda $}$\mathbf{)}$ by (\ref{bc})
because $K(q)$ is invertible in $\mathbf{A}(q)$. Moreover, $v_{\mathbf{%
\lambda }}^{\prime }+qL\mathbf{(}$\textbf{$\lambda $}$\mathbf{)=}v_{\mathbf{%
\lambda }}+qL\mathbf{(}$\textbf{$\lambda $}$\mathbf{)}$ for $v_{\mathbf{%
\lambda }}^{\prime }=v_{\mathbf{\lambda }}+(K(q)-1)v_{\mathbf{\lambda }}$
with $K(q)-1\in qL\mathbf{(}$\textbf{$\lambda $}$\mathbf{)}$. Thus $%
B^{\prime }\mathbf{(}$\textbf{$\lambda $}$\mathbf{)}=B\mathbf{(}$\textbf{$%
\lambda $}$\mathbf{)}$.
\end{proof}

\bigskip

\noindent For any positive integer $n,$ set $V_{n}(\lambda )=U_{q}(\frak{g}%
_{n})\cdot v_{\mathbf{\lambda }}.$ Then $V_{n}(\lambda )$ is a finite
dimensional $U_{q}(\frak{g}_{n})$-module of highest weight $\pi _{n}(\mathbf{%
\lambda })$ and highest vector $v_{\mathbf{\lambda }}$. We denote by $%
(L_{n}(\lambda ),B_{n}(\lambda ))$ the crystal basis of $V_{n}(\lambda )$
obtained as in (\ref{bc}) by authorizing only indices $i_{k}\in I_{n}$.
Clearly, the previous lemma also holds for the finite dimensional highest
weight $U_{q}(\frak{g}_{n})$-modules $V_{n}(\lambda )$. Since $L_{n}(\lambda
)=L\mathbf{(\lambda )}\cap V_{n}(\lambda ),$ the map 
\begin{equation*}
\kappa _{n}:\left\{ 
\begin{array}{l}
B_{n}(\lambda )\rightarrow B\mathbf{(\lambda )} \\ 
\widetilde{f}_{i_{1}}\cdot \cdot \cdot \widetilde{f}_{i_{r}}v_{\mathbf{%
\lambda }}+qL_{n}(\lambda )\mapsto \widetilde{f}_{i_{1}}\cdot \cdot \cdot 
\widetilde{f}_{i_{r}}v_{\mathbf{\lambda }}+qL(\mathbf{\lambda })
\end{array}
\right. 
\end{equation*}
is an embedding of crystals.\ The crystal $B_{n}(\lambda )$ is then
isomorphic to the subcrystal of $B\mathbf{(\lambda )}$ obtained by
considering the vertices which are connected to the highest weight vertex $%
b_{\mathbf{\lambda }}=v_{\mathbf{\lambda }}+qL(\mathbf{\lambda })$ by a path
colored with colors $i\in I_{n}$. In \cite{KN}, Kashiwara and Nakashima have
obtained a natural labelling of the crystal $B_{n}(\lambda )$ by tableaux of
shape $\pi _{n}(\mathbf{\lambda })$. The definition of these tableaux
depends on the root system considered (see Section \ref{sec_plac}).\ For $%
A_{n}$ they coincide with the well known semi-standard tableaux.\ This
labelling of $B_{n}(\lambda )$ induces a natural labelling of $B(\mathbf{%
\lambda })$.\ Indeed consider a vertex $b\in B(\mathbf{\lambda })$ and let $m
$ be the minimal integer such that $b$ belongs to the image of $\kappa _{m}$%
. Then one labels the vertex $b$ by the unique tableau $T$ such that $\kappa
_{m}(T)=b$. The weight of $b$ is then equal to $\mathrm{wt}%
(b)=(d_{1},...,d_{m},\lambda _{m+1},\lambda _{m+2},...)$ where for any $i\in
\{1,...,m\},$ $d_{i}$ is the number of letters $\overline{i}$ minus the
number of letters $i$ in $T$. This means that the vertex labelled by $T$
must in fact be thought as the infinite tableau obtained by adding on the
top of $T,$ first a row of length $\lambda _{m+1}$ containing letters $%
\overline{m+1},$ next a row of length $\lambda _{m+2}$ containing letters $%
\overline{m+2},$ and so on.\ For any positive integer $n,$ the crystal $%
B_{n}(\lambda )$ is obtained by deleting simultaneously in $B(\mathbf{%
\lambda })$ all the arrows colored by $i\geq n$ together with the vertices
they connect. Note that each vertex $b\in B_{n}(\lambda )\subset B(\mathbf{%
\lambda })$ is such that $\widetilde{f}_{n}(b)\neq 0$.

\begin{gather*}
\begin{tabular}{lllllllllllll}
&  &  &  &  &  & $\ \overline{3}$ &  &  &  &  &  &  \\ 
&  &  &  &  &  & {\tiny 2}$\downarrow$ &  &  &  &  &  &  \\ 
&  &  &  &  &  & $\ \overline{2}$ &  &  &  &  &  &  \\ 
&  &  &  &  & $\overset{\text{1}}{\swarrow}$ &  & $\overset{\text{3}}{%
\searrow}$ &  &  &  &  &  \\ 
&  &  &  & $\ \overline{1}$ &  &  &  & 
\begin{tabular}{l}
$\overline{3}$ \\ 
$\overline{2}$%
\end{tabular}
&  &  &  &  \\ 
&  &  &  & {\tiny 0}$\downarrow$ & $\overset{\text{3}}{\searrow}$ &  & $%
\overset{\text{1}}{\swarrow}$ & \ {\tiny 4}$\downarrow$ &  &  &  &  \\ 
&  &  &  & $\ 1$ &  & 
\begin{tabular}{l}
$\overline{3}$ \\ 
$\overline{1}$%
\end{tabular}
&  & 
\begin{tabular}{l}
$\overline{4}$ \\ 
$\overline{3}$ \\ 
$\overline{2}$%
\end{tabular}
&  &  &  &  \\ 
&  &  & $\overset{\text{1}}{\swarrow}$ & {\tiny 3}$\downarrow$ & $\overset{%
\text{0}}{\swarrow}$ & {\tiny 2}$\downarrow$ & $\overset{\text{4}}{\searrow}$
& \ {\tiny 1}$\downarrow$ & $\overset{\text{5}}{\searrow}$ &  &  &  \\ 
&  & $2$ &  & 
\begin{tabular}{l}
$\overline{3}$ \\ 
$1$%
\end{tabular}
&  & 
\begin{tabular}{l}
$\overline{2}$ \\ 
$\overline{1}$%
\end{tabular}
&  & 
\begin{tabular}{l}
$\overline{4}$ \\ 
$\overline{3}$ \\ 
$\overline{1}$%
\end{tabular}
&  & 
\begin{tabular}{l}
$\overline{5}$ \\ 
$\overline{4}$ \\ 
$\overline{3}$ \\ 
$\overline{2}$%
\end{tabular}
&  &  \\ 
$\cdot\cdot\cdot$ &  & $\cdot\cdot\cdot$ &  & $\cdot\cdot\cdot$ &  & $%
\cdot\cdot\cdot$ &  & $\cdot\cdot\cdot$ &  & $\cdot\cdot\cdot$ &  & $%
\cdot\cdot\cdot$%
\end{tabular}
\\
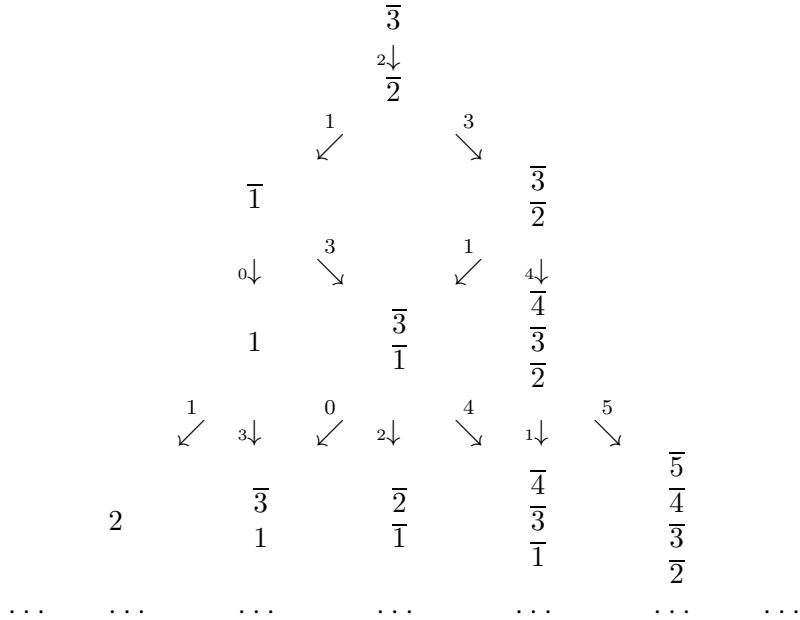
\text{The crystal graph }B^{C_{\infty}}(\mathbf{\omega}_{2})
\end{gather*}
Note that our convention for labelling the crystal graph of the vector
representations are not those used by Kashiwara and Nakashima \cite{KN}.\ To
obtain the original description of $B_{n}(\lambda)$ from that used in the
sequel it suffices to change each letter $k\in\{1,...,n\}$ into $\overline
{n-k+1}$ and each letter $\overline{k}\in\{\overline{1},...,\overline{n}\}$
into $n-k+1.$

\noindent Consider $\mathbf{\lambda}$ and $\mathbf{\mu}$ two dominant
weights of $P^{+}.$ For any positive integer $n$ set 
\begin{multline*}
E_{n}(\mathbf{\lambda},\mathbf{\mu})=\{\mathbf{\nu}\in P^{+}\mid\mathbf{\nu }%
=(\nu_{1},...,\nu_{n},\lambda_{n+1}+\mu_{n+1},\lambda_{n+2}+\mu _{n+2},...)%
\text{ } \\
\text{with }\pi_{n}(\mathbf{\nu})\in P_{n}^{+}\text{ and }\nu_{n}\neq
\lambda_{n}+\mu_{n}\}.
\end{multline*}

\begin{proposition}
With the above notation we have 
\begin{equation*}
\mathcal{V}(\mathbf{\lambda })\otimes \mathcal{V}(\mathbf{\mu })\simeq
\bigoplus_{n\geq 1}\bigoplus_{\mathbf{\nu }\in E_{n}(\mathbf{\lambda },%
\mathbf{\mu })}\mathcal{V}(\mathbf{\nu })^{\oplus m_{\mathbf{\lambda },%
\mathbf{\mu }}^{n,\mathbf{\nu }}}
\end{equation*}
where $m_{\mathbf{\lambda },\mathbf{\mu }}^{n,\mathbf{\nu }}$ is the
multiplicity of $V_{n}(\nu )$ in the tensor product $V_{n}(\lambda )\otimes
V_{n}(\mu )$.
\end{proposition}

\begin{proof}
Suppose that $b=b_{1}\otimes b_{2}\in B(\mathbf{\lambda })\otimes B(\mathbf{%
\mu })$ is of highest weight $\mathbf{\nu }$.\ Then we must have $\widetilde{%
e}_{i}(b)=0$ for any $i\in I.\;$Thus it follows by Lemma \ref{lem_plu_hp}
that $b_{1}$ the highest weight vertex of $B(\mathbf{\lambda })$ (thus has
weight $\mathbf{\lambda }$) and $\varepsilon _{i}(b_{2})\leq \lambda
_{i+1}-\lambda _{i}$ for any $i\in I$.\ Since $b_{2}\in B(\mathbf{\mu })$
there exists a nonnegative integer $n$ such that $\mathrm{wt}(b_{2})=(\delta
_{1},...,\delta _{n},\mu _{n+1},\mu _{n+2},...)$ and $\delta _{n}\neq \mu
_{n}.\;$Then we must have 
\begin{equation*}
\nu =\mathrm{wt}(b_{1})+\mathrm{wt}(b_{2})=(\lambda _{1}+\delta
_{1},...,\lambda _{n}+\delta _{n},\lambda _{n+1}+\mu _{n+1},\lambda
_{n+2}+\mu _{n+2},...).
\end{equation*}
For any $i\in \{1,...,n-1\}$ we have $\left( \lambda _{i+1}+\delta
_{i+1}\right) -\left( \lambda _{i}+\delta _{i}\right) \geq \varepsilon
_{i}(b_{2})+\delta _{i+1}-\delta _{i}.\;$Since $\delta _{i+1}-\delta
_{i}=\varphi _{i}(b_{2})-\varepsilon _{i}(b_{2}),$ we obtain $\ \left(
\lambda _{i+1}+\delta _{i+1}\right) -\left( \lambda _{i}+\delta _{i}\right)
\geq 0$. Thus $\mathbf{\nu }\in E_{n}(\mathbf{\lambda },\mathbf{\mu })$ and $%
\kappa _{n}^{-1}(b)$ is a highest weight vertex of weight $\pi _{n}(\mathbf{%
\lambda })$ in $V_{n}(\lambda )\otimes V_{n}(\mu )$. This means that $%
V_{n}(\lambda )\otimes V_{n}(\mu )$ contains an irreducible component
isomorphic to $V_{n}(\nu )$.

\noindent Conversely, consider $\mathbf{\nu }\in E_{n}(\mathbf{\lambda },%
\mathbf{\mu })$ such that $V_{n}(\lambda )\otimes V_{n}(\mu )$ contains an
irreducible component isomorphic to $V_{n}(\nu ).\;$Let $b^{(n)}$ be the
highest weight vertex of $V_{n}(\lambda )\otimes V_{n}(\mu )$ associated to
this component.\ Then $b=\kappa _{n}(b^{(n)})$ is a highest weight vertex of
weight $\mathbf{\nu }$ in $B(\mathbf{\lambda })\otimes B(\mathbf{\mu }).$
\end{proof}

\begin{example}
The square of the spin representation $\mathcal{V}(\mathbf{\omega }_{0})$ of 
$U_{q}(\frak{g}_{B_{\infty }})$ is isomorphic to the direct sum 
\begin{equation*}
\mathcal{V}(\mathbf{\omega }_{0})^{\otimes 2}\simeq \mathcal{V}(2\mathbf{%
\omega }_{0})\oplus \bigoplus_{n\geq 1}\mathcal{V}(\mathbf{\omega }_{n})
\end{equation*}
\end{example}

\noindent \textbf{Remark: }As in the proof below\textbf{, }the decomposition
of $\mathcal{V}(\mathbf{\lambda })\otimes \mathcal{V}(\mathbf{\mu })$ can be
obtained combinatorially by considering the highest weight vertices of $B(%
\mathbf{\lambda })\otimes B(\mathbf{\mu })$.

\section{The irreducible modules $V_{\infty }(\protect\lambda )$\label%
{sec_modu}}

\subsection{Definition of $V_{\infty}(\protect\lambda)$}

Let $\mathcal{P}_{m}$ denote the set of partitions of length $m\in\mathbb{N}$
and set $\mathcal{P}=\cup_{m\in\mathbb{N}}\mathcal{P}_{m}$. In this paper it
will be convenient to consider the partition $\lambda=(\lambda_{1},...,%
\lambda_{m})\in\mathcal{P}_{m}$ as the \textit{increasing} sequence of
nonnegative integers $0\leq\lambda_{1}\leq\lambda_{2}\leq\cdot\cdot\cdot
\leq\lambda_{m}$.\ For any integer $n\geq m,$ one associates to $\lambda$
the dominant weight 
\begin{equation*}
\phi_{n}(\lambda)=\sum_{i=1}^{m}\lambda_{i}\varepsilon_{i+n-m}\in P_{n}^{+}.
\end{equation*}
Then we write $V_{n}(\lambda)$ for the finite dimensional $U_{q}(\frak{g}%
_{n})$-module of highest weight $\phi_{n}(\lambda).\;$Consider two integers $%
r,s$ such that $m\leq r\leq s.$ The restriction of $V_{s}(\lambda)$ from $%
U_{q}(\frak{g}_{s})$ to $U_{q}(\frak{g}_{r})$ contains a unique irreducible
component isomorphic to $V_{r}(\lambda)$. Thus, there exists a unique
embedding of $U_{q}(\frak{g}_{r})$-modules $f_{s,r}:V_{s}(\lambda)%
\rightarrow V_{r}(\lambda).$ Then we define $V_{\infty}(\lambda)$ as the
direct limit of the direct system $(V_{r}(\lambda),f_{s,r})$ 
\begin{equation*}
V_{\infty}(\lambda)=\lim_{\rightarrow}V_{n}(\lambda).
\end{equation*}
We can suppose without loss of generality that the $U_{q}(\frak{g}_{n})$%
-modules $V_{n}(\lambda)$ are chosen to verify $V_{n}(\lambda)\subset
V_{n+1}(\lambda)$ for $n\geq m.\;$Under this hypothesis, the embeddings $%
f_{s,r}$ define above are trivial and we have $V_{\infty}(\lambda
)=\cup_{n\geq m}V_{n}(\lambda).$

\begin{proposition}
\label{prop_irre}$V_{\infty}(\lambda)$ has the structure of an infinite
dimensional irreducible $U_{q}(\frak{g}_{\infty})$-module.
\end{proposition}

\begin{proof}
Consider a Chevalley generator $g\in\{e_{i},f_{i},t_{i},$ $i\in I\}$ of $%
U_{q}(\frak{g}_{\infty})$ and a vector $v\in V_{\infty}(\lambda)$.\ Then
there exists $N$ such that $v\in V_{N}(\lambda)$ and $g\in U_{q}(\frak{g}%
_{N}).$ The action of $g$ on $v$ is then given by the structure of $U_{q}(%
\frak{g}_{N})$-module on $V_{N}(\lambda).$ Since $V_{n}(\lambda)\subset
V_{n+1}(\lambda)$ for any $n\geq m,$ the action of $g$ on $v$ does not
depend on the integer $N$ considered. One verifies easily that this defines
the structure of a $U_{q}(\frak{g}_{\infty})$-module on $V_{\infty}(\lambda)$%
.

\noindent Suppose that $M$ is a sub-$U_{q}(\frak{g}_{\infty})$-module of $%
V_{\infty}(\lambda).\;$For any positive integer $n\geq m$, write $M_{n}$ for
the restriction of $M\cap V_{n}(\lambda)$ from $U_{q}(\frak{g}_{\infty})$ to 
$U_{q}(\frak{g}_{n})$. Then $M_{n}$ is a $U_{q}(\frak{g}_{n})$-submodule of $%
V_{n}(\lambda)$. Since $V_{n}(\lambda)$ is irreducible, we must have $%
M_{n}=\{0\}$ or $M_{n}=V_{n}(\lambda)$. If $M_{n}=\{0\}$ for each $n\geq m,$
we obtain $M=\{0\}$ because $M=\cup_{n\geq1}M_{n}.\;$Otherwise there exists
an integer $n\geq m$ such that $M_{n}=V_{n}(\lambda)$.$\;$Since $%
M_{n}\subset M_{n+1}$, this yields the equality $M_{p}=V_{p}(\lambda)$ for
any integer $p\geq n$.\ So we obtain $M=U_{n\geq m}M_{n}=V_{\infty}(\lambda)$%
.
\end{proof}

\bigskip

\noindent\textbf{Remarks:}

\noindent$\mathrm{(i):}$ Consider $(V_{r}(\lambda),f_{s,r})$ and $%
(V_{r}^{\prime}(\lambda),f_{s,r}^{\prime})$ two direct systems defined as
above and their direct limits $V_{\infty}(\lambda)=\underrightarrow{\lim}%
V_{n}(\lambda),$ $V_{\infty}^{\prime}(\lambda)=\underrightarrow{\lim}%
V_{n}^{\prime}(\lambda)$. For any $n\geq m,$ let $\Psi_{n}$ be an
isomorphism between $V_{n}(\lambda)$ and $V_{n}^{\prime}(\lambda).$ We must
have $\Psi_{n+1}(V_{n}(\lambda))=V_{n}^{\prime}(\lambda)$ because $%
\Psi_{n+1}(V_{n}(\lambda))$ is a $U_{q}(\frak{g}_{n})$-module isomorphic to $%
V_{n}(\lambda)$ contained in the restriction of $V_{n+1}(\lambda)$ to $U_{q}(%
\frak{g}_{n})$. Thus the map $\Psi:V_{\infty}(\lambda)\rightarrow
V_{\infty}^{\prime}(\lambda)$ such that $\Psi(v)=\Psi_{n}(v)$ for any $v\in
V_{n}(\lambda)$ is a well defined isomorphism of $U_{q}(\frak{g}_{\infty})$%
-modules.

\noindent $\mathrm{(ii):}$ Consider a weight vector $v\in V_{\infty
}(\lambda )$ such that $v\in V_{m}(\lambda )$. Write $\mathrm{wt}%
_{m}(v)=(\delta _{1},...,\delta _{m})\in P_{m}$ for its weight.\ Then for
any $n\geq m,$ $v\in V_{n}(\lambda )$ is a vector of weight $\mathrm{wt}%
_{n}(v)=(\delta _{1},...,\delta _{n},0,...,0)\in P_{n}.$ As a weight vector
of $V_{\infty }(\lambda ),$ the weight of $v$ can be written $\mathrm{wt}%
_{n}(v)=(\delta _{1},...,\delta _{n},0,0,....)\in P$.

\noindent $\mathrm{(iii):}$ For any partition $\lambda $, the weights of $%
V_{\infty }(\lambda )$ are elements of $P$ with a finite nomber of nonzero
coordinates. They should not be confused with the highest weight modules $%
\mathcal{V}(\mathbf{\lambda })$ where $\mathbf{\lambda }\in P^{+}$
introduced in Section \ref{secBa} whose weights have an infinite number of
nonzero coordinates.

\subsection{Littelwood-Richardson rule\label{subsec_lit}}

Consider the infinite ordered alphabets 
\begin{align*}
\mathcal{X}_{A_{\infty}} & =\{\cdot\cdot\cdot<\overline{2}<\overline {1}\},%
\text{\ }\mathcal{X}_{B_{\infty}}=\{\cdot\cdot\cdot<\overline {2}<\overline{1%
}<0<1<2<\cdot\cdot\cdot\} \\
\mathcal{X}_{C_{\infty}} & =\{\cdot\cdot\cdot<\overline{2}<\overline
{1}<1<2<\cdot\cdot\cdot\},\text{\ }\mathcal{X}_{D_{\infty}}=\{\cdot\cdot
\cdot<\overline{2}< 
\begin{array}{l}
1 \\ 
\overline{1}
\end{array}
<2<\cdot\cdot\cdot\}.
\end{align*}
Note that $\mathcal{X}_{D_{\infty}}$ is only partially ordered: $1$ and $%
\overline{1}$ are not comparable. For any positive integer $n,$ we denote by 
$\mathcal{X}_{n}$ the ordered alphabet obtained by considering only the
letters of $\mathcal{X}_{\infty}$ which belong to $\{\overline{n},...,%
\overline{2},\overline{1},0,1,2,...,n\}.$ Set 
\begin{equation*}
\mathrm{wt}(\overline{i})=(\underset{i-1\text{ times}}{\underleftrightarrow
{0,...,0}},1,0,...)\in P,\mathrm{wt}(i)=(\underset{i-1\text{ times}}{%
\underleftrightarrow{0,...,0}},-1,0,...)\in P\text{ and }\mathrm{wt}(0)=0\in
P.
\end{equation*}
The vector representation $V_{\infty}(1)$ of $U_{q}(\frak{g}_{\infty})$ is
the vector space with basis $\{v_{x},$ $x\in\mathcal{X}_{\infty}\}$ where 
\begin{equation*}
q^{h}(v_{x})=q^{<\mathrm{wt}(x),h>}\text{ for }h\in P^{\ast},
\end{equation*}
\begin{equation}
\left\{ 
\begin{tabular}{l}
$f_{i}(v_{\overline{i+1}})=v_{\overline{i}}$, $f_{i}(v_{i})=v_{i+1}$ and $%
f_{i}(v_{x})=0$ if $x\notin\{\overline{i+1},i\}$ \\ 
$e_{i}(v_{\overline{i}})=v_{\overline{i+1}}$, $f_{i}(v_{i+1})=v_{i}$ and $%
e_{i}(v_{x})=0$ if $x\notin\{\overline{i},i+1\}$%
\end{tabular}
\right. \text{ if }i\neq0  \label{action fi_vect}
\end{equation}
and the action of $e_{0},f_{0}$ depends on the type considered:

\begin{align}
& \left\{ 
\begin{tabular}{l}
$e_{0}(v_{1})=v_{0},$ $e_{0}(v_{0})=(q+q^{-1})v_{\overline{1}}$ and $%
e_{0}(v_{x})=0$ if $x\notin\{0,1\}$ \\ 
$f_{0}(v_{\overline{1}})=v_{0},\text{ }f_{0}(v_{0})=(q+q^{-1})v_{1}\text{
and }f_{0}(v_{x})=0\text{ if }x\notin\{\overline{1},0\}$%
\end{tabular}
\right. \text{ for type }B_{\infty}  \label{actionf0} \\
& \left\{ 
\begin{tabular}{l}
$e_{0}(v_{1})=v_{\overline{1}}$ and $e_{0}(v_{x})=0$ if $x\neq1$ \\ 
$f_{0}(v_{\overline{1}})=v_{1}\text{ and }f_{0}(v_{x})=0\text{ if }x\neq%
\overline{1}$%
\end{tabular}
\right. \text{ for type }C_{\infty}  \notag \\
& \left\{ 
\begin{tabular}{l}
$f_{0}(v_{\overline{2}})=v_{1}$, $f_{0}(v_{\overline{i}})=v_{2}$ and $%
f_{0}(v_{x})=0$ if $x\notin\{\overline{2},\overline{1}\}$ \\ 
$e_{i}(v_{2})=v_{\overline{1}}$, $f_{i}(v_{1})=v_{\overline{2}}$ and $%
e_{0}(v_{x})=0$ if $x\notin\{1,2\}$%
\end{tabular}
\right. \text{ for type }D_{\infty}  \notag
\end{align}
For any positive integer $n,$ the $U_{q}(\frak{g}_{n})$-module $V_{n}(1)$ is
identified with the vector space with basis $\{v_{x},$ $x\in\mathcal{X}%
_{n}\} $ with the action of the Chevalley generators $e_{i},f_{i},t_{i},$ $%
i\in I_{n}$.

\bigskip

\noindent For any nonnegative integer $\ell,$ there exists a finite number
of irreducible modules $V_{n}^{(1)},...,V_{n}^{(k)}$ such that 
\begin{equation}
V_{n}(1)^{\otimes\ell}=\bigoplus_{k=1}^{r}V_{n}^{(k)}.  \label{deco}
\end{equation}
Set $E_{n,\ell}=\{k\in\{1,...,r\}\mid V_{n}^{(k)}\simeq V_{n}(\lambda)$ with 
$\left| \lambda\right| =\ell\}$ and $F_{n,\ell}=\{1,...,r\}-E_{n,\ell}.$
Then consider the sub-modules $M_{n,\ell}$ and $N_{n,\ell}$ defined by 
\begin{equation*}
M_{n,\ell}=\bigoplus_{k\in E_{n,\ell}}V_{n}^{(k)}\text{ and }N_{n,\ell
}=\bigoplus_{k\in F_{n,\ell}}V_{n}^{(k)}.
\end{equation*}
These modules do not depend on the decomposition (\ref{deco}) and we have $%
V_{n}(1)^{\otimes\ell}=M_{n,\ell}\oplus N_{n,\ell}.$ The module $M_{n,\ell}$
(resp.\ $N_{n,\ell})$ contains all the highest weight vectors of highest
weight $\lambda$ verifying $\left| \lambda\right| =\ell$ (resp. $\left|
\lambda\right| <\ell$). Note that $N_{n,\ell}=\{0\}$ only for type $A_{n-1}$.

\noindent The crystal basis of $V_{n}(1)$ is the pair $(L_{n}(1),B_{n,1})$
where 
\begin{equation*}
L_{n}(1)=\bigoplus_{x\in\mathcal{X}_{n}}\mathbf{A}(q)v_{x}
\end{equation*}
and for types $A_{n-1},B_{n},C_{n},D_{n}$, the crystals $B_{n,1}$ are
respectively 
\begin{gather*}
\overline{n}\overset{n-1}{\rightarrow}\overline{n-1}\overset{n-2}{%
\rightarrow }\cdot\cdot\cdot\cdot\rightarrow\overline{2}\overset{1}{%
\rightarrow}\overline{1} \\
\overline{n}\overset{n-1}{\rightarrow}\overline{n-1}\overset{n-2}{%
\rightarrow }\cdot\cdot\cdot\cdot\rightarrow\overline{2}\overset{1}{%
\rightarrow}\overline{1}\overset{0}{\rightarrow}0\overset{0}{\rightarrow}1%
\overset{1}{\rightarrow}2\cdot\cdot\cdot\cdot\overset{n-2}{\rightarrow}n-1%
\overset{n-1}{\rightarrow}n \\
\overline{n}\overset{n-1}{\rightarrow}\overline{n-1}\overset{n-2}{%
\rightarrow }\cdot\cdot\cdot\cdot\rightarrow\overline{2}\overset{1}{%
\rightarrow}\overline{1}\overset{0}{\rightarrow}1\overset{1}{\rightarrow}%
2\cdot\cdot \cdot\cdot\overset{n-2}{\rightarrow}n-1\overset{n-1}{\rightarrow}%
n \\
\overline{n}\overset{n-1}{\rightarrow}\overline{n-1}\overset{n-2}{%
\rightarrow }\cdot\cdot\cdot\overset{3}{\rightarrow}\overline{3}\overset{2}{%
\rightarrow} 
\begin{tabular}{c}
$1$ \ \  \\ 
\ \ $\overset{0}{\nearrow}$ $\ \ \ \overset{\text{ \ \ \ }1}{\text{ \ }%
\searrow}$ \ \ \  \\ 
$\overline{2}\ \ \ \ \ \ \ \ \ \ \ \ \ \ \ \ \ \ 2$ \\ 
\ $\underset{1\text{ \ \ \ }}{\searrow}$ \ \ \ $\underset{0}{\nearrow}$ \ \
\  \\ 
$\overline{1}$ \ 
\end{tabular}
\overset{2}{\rightarrow}3\overset{3}{\rightarrow}\cdot\cdot\cdot\overset{n-2%
}{\rightarrow}n-1\overset{n-1}{\rightarrow}n.
\end{gather*}
The crystal basis of $V_{n}(1)^{\otimes\ell}$ can be realized as the pair $%
(L_{n,\ell},B_{n,\ell})$ with $L_{n,\ell}=\oplus\mathbf{A}%
(q)v_{x_{1}}\otimes\cdot\cdot\cdot\otimes v_{x_{\ell}}$ and $%
B_{n,\ell}=B_{n,1}^{\otimes\ell}$. The weight of the vertex $%
b=x_{1}\otimes\cdot\cdot \cdot\otimes x_{\ell}$ is the $n$-tuple $\mathrm{wt}%
(b)=(d_{1},...,d_{n})\in P_{n}$ where for any $i\in I_{n},$ $d_{i}$ is the
number of letters $x_{k}=\overline{i}$ in $b$ minus its number of letters $%
x_{k}=i$. In particular $\left| \mathrm{wt}(b)\right| $ is the number of
barred letters in $b$ minus its number of unbarred letters. The pair 
\begin{equation}
L_{n,\ell}\cap N_{n,\ell}\text{, }B_{n,l}^{N}=\{b+qL_{n,\ell}\cap N_{n,\ell
}\mid b+qL_{n,\ell}\in B_{n,\ell},b\in L_{n,\ell}\cap N_{n,\ell}\}
\label{cbL}
\end{equation}
is a crystal basis of $N_{n,\ell}$. This implies that $B_{n,\ell}^{N}$ is
the subcrystal of $B_{n,\ell}$ containing all the connected components
isomorphic to $B_{n}(\lambda)$ with $\left| \lambda\right| <\ell$.

\noindent Given $m\leq n$ two integers, denote by $\rho_{m,n}:V_{n}(1)^{%
\otimes\ell}\rightarrow V_{m}(1)^{\otimes\ell}$ the canonical projection.
Then the image of the crystal basis $(L_{n,\ell},B_{n,\ell})$ by $\rho_{m,n}$
is the crystal basis $(L_{m,\ell},B_{m,\ell})$.

\begin{lemma}
Consider $m$ a nonnegative integer. Then for $n\geq m$ sufficiently large $%
B_{m,\ell}\cap B_{n,\ell}^{N}=\emptyset$.
\end{lemma}

\begin{proof}
For any $b\in B_{m,\ell}$, denote by $B_{n,\ell}(b)$ the connected component
of $B_{n,\ell}$ containing $b$. We are going to prove by induction on $\ell$
that for $n$ sufficiently large, there exists a partition $\lambda$ such
that $\mathrm{wt}(b_{h})=\phi_{n}(\lambda)$ with $\left| \lambda\right|
=\ell $.

\noindent For $\ell=1,$ this is immediate because $b_{h}=\overline{n}$ has
weight $\phi_{n}(1).$ So, consider $b\in B_{m,\ell},\ell\geq2$ and write $%
b=b^{\flat}\otimes x_{\ell},b_{h}=b_{h}^{\flat}\otimes x_{\ell,h}$ where $%
b^{\flat},b_{h}^{\flat}\in B_{n,\ell-1}$ and $x_{\ell},x_{\ell,h}\in B_{n,1}$%
. Then it follows from (\ref{TENS2}) and Lemma \ref{lem_plu_hp} that $%
b_{h}^{\flat}$ is the highest weight vertex of $B_{n,\ell-1}(b^{\flat})$.
Choose $n\geq\ell$ sufficiently large so that there exists a partition $%
\gamma$ verifying $\mathrm{wt}(b^{\flat})=\phi_{n}(\gamma)$ with $\left|
\gamma\right| =\ell-1$. The letters of $b_{h}^{\flat}$ belong to $\{%
\overline{n},...,\overline{n-\ell+2}\}$ for $\left| \gamma\right| =\ell-1$
and $\phi_{n}(\gamma)\in P_{n-1}^{+}$. Moreover for any $i\in I_{n},$ we
must have $\varepsilon_{i}(x_{\ell,h})\leq\varphi_{i}(b_{h}^{\flat})$. Since 
$b\in B_{m,\ell}$ and $x_{\ell,h}\leq x_{\ell},$ we have $x_{\ell,h}\leq m$.
Suppose that $x_{\ell,h}=k+1$ with $k\in\{0,...,m-1\}$. Then $%
\varepsilon_{k}(x_{\ell,h})=1,$ thus $\varphi_{k}(b_{h}^{\flat})\geq1$. This
implies that $b_{h}^{\flat}$ contains at least a letter $\overline{k+1}$.
Hence we obtain $n-\ell+1\leq k\leq m-1$.\ This means that $n$ can be chosen
sufficiently large so that $b_{h}^{\flat}$ has weight $\phi_{n}(\gamma)$
with $\gamma\in\mathcal{P}$, $\left| \gamma\right| =\ell-1$ and $x_{\ell,h}$
is barred. For such an integer $n,$ $b_{h}$ contains only barred letters.
Thus there exists a partition $\lambda$ such that $\phi_{n}(\lambda)=\mathrm{%
wt}(b_{h})$ and $\left| \lambda\right| =\ell$.
\end{proof}

\begin{lemma}
\label{lem_fund}For any positive integer $m,$ there exists an integer $n\geq
m$ such that $V_{m}(1)^{\otimes\ell}$ is a submodule of the restriction of $%
M_{n,\ell}$ to $U_{q}(\frak{g}_{m}).$
\end{lemma}

\begin{proof}
From the decomposition $V_{n}(1)^{\otimes\ell}=M_{n,\ell}\oplus N_{n,\ell},$
we obtain $V_{m}(1)^{\otimes\ell}=M_{n,\ell}\cap
V_{m}(1)^{\otimes\ell}\oplus N_{n,\ell}\cap V_{m}(1)^{\otimes\ell}$ for any $%
m\leq n$. Since $L_{m,\ell }=L_{n,\ell}\cap V_{m}(1)^{\otimes\ell}$, we
derive from (\ref{cbL}) that the pair 
\begin{equation*}
L_{m,\ell}\cap N_{n,\ell}\text{ and }B_{m,n,l}^{N}=\{b+qL_{m,\ell}\cap
N_{n,\ell}\mid b+qL_{n,\ell}\cap N_{n,\ell}\in B_{n,\ell}^{N},b\in L_{m,\ell
}\cap N_{n,\ell}\}
\end{equation*}
is a crystal basis of the $U_{q}(\frak{g}_{m})$-module $N_{n,\ell}\cap
V_{m}(1)^{\otimes\ell}$.\ In particular, the vertices of $B_{m,n,l}^{N}$ are
the vertices of $B_{m,\ell}$ which belongs to $B_{n,\ell}^{N}$. By the above
Lemma, we obtain that $B_{m,n,l}^{N}=\emptyset$ for $n$ sufficiently large.
For such an integer $n,$ we must have $N_{n,\ell}\cap V_{m}(1)^{\otimes\ell
}=\{0\}$ because the dimension of $N_{n,\ell}\cap V_{m}(1)^{\otimes\ell}$ is
equal to the number of vertices of $B_{m,n,l}^{N}$. Finally, this gives $%
V_{m}(1)^{\otimes\ell}=M_{n,\ell}\cap V_{m}(1)^{\otimes\ell},$ thus $%
V_{m}(1)^{\otimes\ell}$ is a submodule of the restriction of $M_{n,\ell}$ to 
$U_{q}(\frak{g}_{m}).$
\end{proof}

\bigskip

\noindent For any $x\in\mathcal{X}$, set 
\begin{equation*}
\theta(x)=\left\{ 
\begin{array}{l}
i+1\text{ if }x=i\text{ is unbarred} \\ 
\overline{i+1}\text{ if }x=\overline{i}\text{ is barred}
\end{array}
\right. .
\end{equation*}
We also denote by $\theta$ the linear map defined on $V_{\infty}(1)^{\otimes%
\ell}$ by setting $\theta(v_{x_{1}}\otimes\cdot\cdot\cdot\otimes
v_{x_{\ell}})=v_{\theta(x_{1})}\otimes\cdot\cdot\cdot\otimes v_{\theta
(x_{\ell})}$.

\begin{lemma}
\label{lem-trans}Let $v$ be vector in $V_{\infty}(1)^{\otimes\ell}$ of
weight $\delta=(\delta_{1},...,\delta_{n},0,...)\in P$ where $%
\delta_{n}\neq0.$ Suppose that $\left| \delta\right| =\ell$. Then, there
exists $K(q)$ in $\mathbb{Q}(q)$ such that 
\begin{equation*}
f_{n}^{\delta_{n}}\cdot\cdot\cdot f_{1}^{\delta_{1}}(\theta(v))=K(q)v
\end{equation*}
\end{lemma}

\begin{proof}
Observe first that $\theta(v)\in V_{\infty}(1)^{\otimes\ell}$ is a vector of
weight $(0,\delta_{1},...,\delta_{n},0,...)\in P.$ Since $%
\delta_{1}+\cdot\cdot\cdot+\delta_{n}=\ell,$ each tensor $%
v_{x_{1}}\otimes\cdot \cdot\cdot\otimes v_{x_{\ell}}$ appearing in its
decomposition 
\begin{equation}
\theta(v)=\sum a_{x_{1}\cdot\cdot\cdot x_{\ell}}(q)v_{x_{1}}\otimes\cdot
\cdot\cdot\otimes v_{x_{\ell}}  \label{dec_teta}
\end{equation}
on the tensor basis of $V_{\infty}(1)^{\otimes\ell}$ must verify $x_{k}\in\{%
\overline{2},...,\overline{n+1}\}$. Then by (\ref{tensor1}) there exists $%
K_{1}(q)$ in $\mathbb{Q}(q)$ such that 
\begin{equation*}
f_{1}^{\delta_{1}}(\theta(v))=K_{1}(q)v^{(1)}
\end{equation*}
where $v^{(1)}$ is obtained by replacing each factor $v_{\overline{2}}$
appearing in (\ref{dec_teta}) by a factor $v_{\overline{1}}$. Indeed the
factors $v_{x_{k}}$ with $x_{k}\neq\overline{2}$ do not interfere with the
computation of $f_{1}^{\delta_{1}}(v)$. Similarly the tensors appearing in
the decomposition of $v^{(1)}$ on the tensor basis have factors in $%
\{\overline {1},\overline{3},...,\overline{n+1}\}.\;$Hence there also exists 
$K_{2}(q)$ in $\mathbb{Q}(q)$ such that$f_{2}^{%
\delta_{2}}(v^{(1)})=K_{2}(q)v^{(2)}$ where $v^{(2)}$ is obtained by
replacing each factor $v_{\overline{3}}$ appearing in $v^{(1)}$ by a factor $%
v_{\overline{2}}$ (the other factors do not interfere with the computation
of $f_{2}^{\delta_{2}}(v^{(1)})$). The lemma follows by an easy induction.
\end{proof}

\begin{theorem}
\label{th_decovec}For any integer $\ell$, $V_{\infty}(1)^{\otimes\ell}$
decomposes into irreducible components following 
\begin{equation}
V_{\infty}(1)^{\otimes\ell}\simeq\bigoplus_{\nu\in\mathcal{P},\left|
\nu\right| =\ell}V_{\infty}(\nu)^{\oplus n_{\nu}}  \label{decvec}
\end{equation}
where $n_{\nu}$ is the number of standard tableaux of shape $\nu$.
\end{theorem}

\begin{proof}
The dimension of the space generated by the highest weight vectors in $%
V_{\ell}(1)^{\otimes\ell}$ of weight $\phi_{\ell}(\nu)\in P_{\ell}$ such
that $\left| \nu\right| =\ell$ is equal to the number of standard tableaux
of shape $\nu$.\ Fix a basis $\{v_{t}^{\ell}\mid t\in\mathrm{ST}(\ell)\}$
(labelled by the standard tableaux with $\ell$ boxes) of highest weight
vectors in $V_{\ell}(1)^{\otimes\ell}$ having a weight of the form $\phi
_{\ell}(\nu)$ with $\left| \nu\right| =\ell$. Consider an integer $n\geq
\ell $ and a standard tableau $t\in\mathrm{ST}(\ell)$ of shape $\nu.\;$Since 
$\left| \nu\right| =\ell,$ the decomposition of $v_{t}^{\ell}$ on the tensor
basis of $V_{\ell}(1)^{\otimes\ell}$ makes only appear tensors $%
v_{x_{1}}\otimes\cdot\cdot\cdot\otimes v_{x_{\ell}}$ such that $x_{k}$ is
barred for any $k\in\{1,...,\ell\}$.

\noindent By (\ref{tensor2}), $\theta^{(n-\ell)}(v_{t}^{\ell})=v_{t}^{n}$ is
a highest weight vector of $V_{n}(1)^{\otimes\ell}$ of weight $%
\phi_{n}(\nu)\in P_{n}$ such that $\left| \nu\right| =\ell$.\ This permits
to define a $U_{q}(\frak{g}_{n})$-module $V_{n}(t)=U_{q}(\frak{g}_{n})\cdot
v_{t}^{n}$ isomorphic to $V_{n}(\nu)$. By Lemma \ref{lem-trans} there exists 
$K(q)\in\mathbb{Q}(q)$ such that 
\begin{equation*}
f_{n}^{\nu_{p}}\cdot\cdot\cdot
f_{n-p+1}^{\nu_{1}}(v_{t}^{n+1})=K(q)v_{t}^{n}.
\end{equation*}
This gives the inclusion $V_{n}(t)\subset V_{n+1}(t).\;$Hence $V_{\infty
}(t)=\cup_{n\geq\ell}V_{n}(t)$ is isomorphic to $V_{\infty}(\nu)$.

\noindent Since $\{v_{t}^{\ell}\mid t\in\mathrm{ST}(\ell)\}$ is a basis of
highest weight vectors in $V_{\ell}(1)^{\otimes\ell}$ of weight $\phi_{\ell
}(\nu)$ with $\left| \nu\right| =\ell$, for any $n\geq\ell,$ $%
\{v_{t}^{n}\mid t\in\mathrm{ST}(\ell)\}$ is a basis of highest weight
vectors in $V_{n}(1)^{\otimes\ell}$ of weight $\phi_{n}(\nu)$ with $\left|
\nu\right| =\ell$. Thus $M_{n,\ell}=\oplus_{t\in\mathrm{ST}(\ell)}V_{n}(t).$
Set $M_{\ell}=\oplus_{t\in\mathrm{ST}(\ell)}V_{\infty}(t)$. Then $%
M_{\ell}=\cup_{n\geq\ell}M_{n,\ell}$. For any integer $m\geq\ell,$ we derive
by Lemma \ref{lem_fund} that there exists an integer $n\geq m$ such that $%
V_{m}(1)^{\otimes\ell}\subset M_{n,\ell}\subset M_{\ell}$.\ Since $V_{\infty
}(1)^{\otimes\ell}=\cup_{m\geq\ell}V_{m}(1)^{\otimes\ell},$ this gives $%
V_{\infty}(1)^{\otimes\ell}\subset M_{\ell}\subset V_{\infty}(1)^{\otimes
\ell}$.\ Hence $V_{\infty}(1)^{\otimes\ell}=M_{\ell}$ and the theorem is
proved for $M_{\ell}$ is isomorphic by definition to $\oplus_{\nu \in%
\mathcal{P},\left| \nu\right| =\ell}V_{\infty}(\nu)^{\oplus n_{\nu}}$.
\end{proof}

\bigskip

\noindent For any $n\geq\ell,$ let $\rho_{n}$ be the canonical projection $%
\rho_{n}:V_{\infty}(1)^{\otimes\ell}\rightarrow V_{n}(1)^{\otimes\ell}.$ We
derive from the above proof the decomposition 
\begin{equation}
V_{n}(1)^{\otimes\ell}=\bigoplus_{t\in\mathrm{ST}(\ell)}\rho_{n}(V_{\infty
}(t))=\bigoplus_{t\in\mathrm{ST}(\ell)}V_{\infty}(t)\cap V_{n}(1)^{\otimes
\ell}.  \label{dec_ron}
\end{equation}
The irreducible components $V_{n}(\lambda)$ appearing in the decomposition
of $V_{n}(1)^{\otimes\ell}$ verify $\left| \lambda\right| \leq\ell.$
Moreover the multiplicity of $V_{n}(\nu)$ with $\left| \nu\right| =\ell$ in $%
V_{n}(1)^{\otimes\ell}$ is equal to $n_{\nu}$.\ Since $v_{t}^{n}$ belongs to 
$V_{\infty}(t)\cap V_{n}(1)^{\otimes\ell}$, the $U_{q}(\frak{g}_{n})$-module 
$\rho_{n}(V_{\infty}(t))=V_{\infty}(t)\cap V_{n}(1)^{\otimes\ell}$ contains
an irreducible component isomorphic to $V_{n}(\nu)$ where $\nu$ is the shape
of $t$. Thus, by (\ref{dec_ron}), the multiplicity of $V_{n}(\nu)$ in $\rho
_{n}(V_{\infty}(t))$ is equal to $1$. We deduce that for any partition $\nu$
such that $\left| \nu\right| =\ell,$ the decomposition of $\rho
_{n}(V_{\infty}(\nu))$ has the form 
\begin{equation}
\rho_{n}(V_{\infty}(t))\simeq V_{n}(\nu)\oplus\bigoplus_{\lambda\in \mathcal{%
P}\left| \lambda\right| <\ell}V_{n}(\lambda)^{\oplus m_{\nu ,\lambda}}.
\label{dec_Vinf}
\end{equation}

\begin{corollary}
\label{dec_tens}Consider $\lambda$ and $\mu$ two partitions respectively of
lengths $m$ and $n.\;$Then we have the following decomposition 
\begin{equation*}
V_{\infty}(\lambda)\otimes V_{\infty}(\mu)\simeq\bigoplus_{\nu\in\mathcal{P}%
}V_{\infty}(\nu)^{\oplus c_{\lambda,\mu}^{\nu}}
\end{equation*}
where $c_{\lambda,\mu}^{\nu}$ is the Littelwood-Richardson coefficient
associated to the partitions $\lambda,\mu$ and $\nu$.
\end{corollary}

\begin{proof}
Set $\ell=\left| \lambda\right| +\left| \nu\right| .$ Since $V_{\infty
}(\lambda)$ and $V_{\infty}(\mu)$ can be realized as irreducible components
respectively of $V_{\infty}(1)^{\otimes\left| \lambda\right| }$ and $%
V_{\infty}(1)^{\otimes\left| \mu\right| }$, we deduce by the previous
theorem that the tensor product $V_{\infty}(\lambda)\otimes
V_{\infty}(\mu)\subset V_{\infty}(1)^{\otimes\ell}$ decomposes in the form 
\begin{equation}
V_{\infty}(\lambda)\otimes V_{\infty}(\mu)\simeq\bigoplus_{\nu\in\mathcal{P}%
}V_{\infty}(\nu)^{\oplus d_{\lambda,\mu}^{\nu}}.  \label{dec_tenslmu}
\end{equation}
For any $\nu\in\mathcal{P}$ such that $\left| \nu\right| =\ell$ and any
integer $n\geq\ell$, the multiplicity of $V_{n}(\nu)$ in the tensor product $%
V_{n}(\lambda)\otimes V_{n}(\mu)$ is the Littlewood-Richardson coefficient $%
c_{\lambda,\mu}^{\nu}$.\ It is well known that $c_{\lambda,\mu}^{n}$ is
equal to the cardinality of the set $\mathbf{LR}_{\lambda,\mu}^{\nu}$ of
Littlewood-Richardson tableaux of shape $\nu/\lambda$ and content $\mu$.
(see \cite{Fu}). Fix a basis $\{v_{t}^{\ell}\mid t\in\mathbf{LR}%
_{\lambda,\mu}^{\nu}\}$ of highest weight vectors of weight $%
\phi_{\ell}(\nu) $ in $V_{\ell }(\lambda)\otimes V_{\ell}(\mu)\subset
V_{\ell}(1)^{\otimes\ell}$. For any $n\geq
l,v_{t}^{n}=\theta^{(n-\ell)}(v_{t}^{\ell})$ is a highest weight vector of
weight $\phi_{n}(\nu)$ in $V_{n}(\lambda)\otimes V_{n}(\mu)\subset
V_{n}(1)^{\otimes\ell}.$ Then $V_{n}(t)=U_{q}(\frak{g}_{n})\cdot v_{t}^{n}$
is isomorphic to $V_{n}(\nu)$. Similarly to the above proof we can define in 
$V_{\infty}(\lambda)\otimes V_{\infty}(\mu)$ irreducible $U_{q}(\frak{g}%
_{\infty})$-submodules $V_{\infty}(t)$ isomorphic to $V_{\infty}(\nu)$ so
that $\oplus_{t\in\mathbf{LR}_{\lambda,\mu}^{\nu}}V_{\infty}(t)$ is a $U_{q}(%
\frak{g}_{\infty})$-submodule of $V_{\infty}(\lambda)\otimes V_{\infty
}(\mu).$ This shows that $d_{\lambda,\mu}^{\nu}\geq c_{\lambda,\mu}^{\nu}$.
For any $\nu\in\mathcal{P}$ such that $\left| \nu\right| =\ell$ and any
integer $n\geq\ell$, we derive from (\ref{dec_Vinf}) that $%
\rho_{n}(V_{\infty }(\lambda)\otimes
V_{\infty}(\mu))=\rho_{n}(V_{\infty}(\lambda))\otimes
\rho_{n}(V_{\infty}(\mu))$ contains $c_{\lambda,\mu}^{\nu}$ irreducible
components isomorphic to $V_{n}(\nu)$. By (\ref{dec_tenslmu}), this implies
that $d_{\lambda,\mu}^{\nu}\leq c_{\lambda,\mu}^{\nu}.$
\end{proof}

\bigskip

\noindent In the sequel we denote by $\mathcal{M}$ the category of modules $%
M $ which can be realized as submodules of a tensor power $%
V_{\infty}(1)^{\otimes\ell}.$ By corollary $\mathcal{M}$ is closed under
taking direct sums or tensor product of finitely many $U_{q}(\frak{g}%
_{\infty})$-modules. Moreover any module in $\mathcal{M}$ is isomorphic to a
direct sum of irreducible $U_{q}(\frak{g}_{\infty})$-modules $%
V_{\infty}(\lambda)$ with $\lambda\in\mathcal{P}.$

\subsection{Explicit realization of the modules $V_{\infty}(\protect\lambda )%
\label{subsec_qwed}$}

\noindent By Theorem \ref{th_decovec}, we have the decomposition 
\begin{equation}
V_{\infty}(1)^{\otimes2}\cong V_{\infty}(2,0)\oplus V_{\infty}(1,1)
\label{dec_square_of B(1)}
\end{equation}
For any positive integer $p,$ set 
\begin{equation*}
N_{p}=\underset{i=0}{\overset{p-2}{\sum}}V_{\infty}(1)^{\otimes i}\otimes
V_{\infty}(2,0)\otimes V_{\infty}(1)^{\otimes(p-2-i)}
\end{equation*}
(where the sum is not direct).

\begin{proposition}
With the above notation 
\begin{equation*}
V_{\infty}(1^{p})\simeq V_{\infty}(1)^{\otimes p}/N_{p}
\end{equation*}
\end{proposition}

\begin{proof}
By Theorem \ref{th_decovec} and Corollary \ref{dec_tens} we have 
\begin{equation*}
V_{\infty}(1)^{\otimes p}\simeq\bigoplus_{\nu\in\mathcal{P}%
_{p}}V_{\infty}(\nu)^{\oplus n_{\nu}}\text{ and }N_{p}\simeq\bigoplus_{\nu\in%
\mathcal{P}_{p}^{\ast}}V_{\infty}(\nu)^{\oplus n_{\nu}}
\end{equation*}
where $\mathcal{P}_{p}^{\ast}$ is the set of partitions of length $p$ with
at least a part $\geq2.\;$Since $\mathcal{P}_{p}=\mathcal{P}%
_{p}^{\ast}\cup\{(1^{p})\},$ this implies the proposition.
\end{proof}

\bigskip

\noindent The representation $V_{\infty}(1^{p})$ can thus be regarded as a $%
q $-analogue of the $p$-th wedge product of $V_{\infty}(1).\;$Let $\Psi_{p}$
be the canonical projection $V_{\infty}(1)^{\otimes p}\rightarrow V_{\infty
}(1^{p})$.\ Denote by $v_{x_{1}}\wedge\cdot\cdot\cdot\wedge v_{x_{p}}$ the
image of the vector $v_{x_{1}}\otimes\cdot\cdot\cdot\otimes v_{x_{p}}$ by $%
\Psi_{p}$. The following lemma has been proved in \cite{lec3} (see also \cite
{JMO}):

\begin{lemma}
\ \ \ \ \ \ \ \ \ \ \label{lem_relation}

\noindent In $V_{B_{\infty}}(1^{2})$ we have the relations

\begin{enumerate}
\item  $v_{x}\wedge v_{x}=0$ for $x\neq0$,

\item  $v_{y}\wedge v_{x}=-q^{2}v_{x}\wedge v_{y}$ for $x\neq\overline{y}$
and $x<y,$

\item  $v_{i}\wedge v\overline{_{i}}=-q^{4}v_{\overline{i}}\wedge
v_{i}+(1-q^{4})\underset{k=1}{\overset{i-1}{\sum}}(-1)^{i-k}q^{2(i-k)}v_{%
\overline{k}}\wedge v_{k}+(-1)^{i}q^{2i-1}v_{0}\wedge v_{0}$ for $i=1,...,n.$
\end{enumerate}

\noindent In $V_{C_{\infty}}(1^{2})$ we have the relations

\begin{enumerate}
\item  $v_{x}\wedge v_{x}=0$ for $x\neq0$,

\item  $v_{y}\wedge v_{x}=-qv_{x}\wedge v_{y}$ for $x\neq\overline{y}$ and $%
x<y,$

\item  $v_{i}\wedge v\overline{_{i}}=-q^{2}v_{\overline{i}}\wedge
v_{i}+(1-q^{2})\underset{k=1}{\overset{i-1}{\sum}}(-1)^{i-k}q^{i-k}v_{%
\overline{k}}\wedge v_{k}$ for $i=1,...,n.$
\end{enumerate}

\noindent In $V_{D_{\infty}}(1^{2})$ we have the relations

\begin{enumerate}
\item  $v_{x}\wedge v_{x}=0$ for $x\neq0$,

\item  $v_{y}\wedge v_{x}=-qv_{x}\wedge v_{y}$ for $x\neq\overline{y}$ and $%
x<y,$

\item  $v_{i}\wedge v\overline{_{i}}=-q^{2}v_{\overline{i}}\wedge
v_{i}+(1-q^{2})\underset{k=2}{\overset{i-1}{\sum}}(-1)^{i-k}q^{i-k}v_{%
\overline{k}}\wedge v_{k}+(-q)^{i-1}(v_{1}\wedge v_{\overline{1}}+v_{%
\overline{1}}\wedge v_{1})$ for $i=2,...,n.$
\end{enumerate}
\end{lemma}

\noindent The columns of types $A_{\infty},B_{\infty},C_{\infty}$ and $%
D_{\infty}$ are Young diagrams $C$ of shape column respectively of the form 
\begin{equation}
C= 
\begin{tabular}{|l|}
\hline
$C_{-}$ \\ \hline
\end{tabular}
,C= 
\begin{tabular}{|l|}
\hline
$C_{-}$ \\ \hline
$C_{0}$ \\ \hline
$C_{+}$ \\ \hline
\end{tabular}
,C= 
\begin{tabular}{|l|}
\hline
$C_{-}$ \\ \hline
$C_{+}$ \\ \hline
\end{tabular}
\text{ and }C= 
\begin{tabular}{|c|}
\hline
$D_{-}$ \\ \hline
$D$ \\ \hline
$D_{+}$ \\ \hline
\end{tabular}
\label{col}
\end{equation}
where $C_{-},C_{+},C_{0},D_{-},D_{+}$ and $D$ are column shaped Young
diagrams such that 
\begin{equation*}
\left\{ 
\begin{tabular}{l}
$C_{-}$ is filled by strictly increasing barred letters from top to bottom
\\ 
$C_{+}$ is filled by strictly increasing unbarred letters from top to bottom
\\ 
$C_{0}$ is filled by letters $0$ \\ 
$D_{-}$ is filled by strictly increasing letters $\leq\overline{2}$ from top
to bottom \\ 
$D_{+}$ is filled by strictly increasing letters $\geq2$ from top to bottom
\\ 
$D$ is filled by letters $\overline{1}$ or $1$ with different letters in two
adjacent boxes
\end{tabular}
\right. .
\end{equation*}
We denote by $\mathbf{C}(p)$ the set of columns of height $p.\;$The reading
of the column $C$ is the word $\mathrm{w}(C)$ obtained by reading the
letters of $C$ from top to bottom. For any column $C$ with $\mathrm{w}%
(C)=c_{1}\cdot \cdot\cdot c_{p}$ where the $c_{i}$'s are letters, we set $%
v_{C}=v_{c_{1}}\wedge\cdot\cdot\cdot\wedge v_{c_{p}}$.\ Then each vector $%
\Psi_{p}(v_{x_{1}}\otimes\cdot\cdot\cdot\otimes
v_{x_{p}})=v_{x_{1}}\wedge\cdot\cdot\cdot\wedge v_{x_{p}}$ can be decomposed
into a linear combination of vectors $v_{C}$ by applying from left to right
a sequence of relations given in the above Lemma.

\begin{lemma}
The vectors of $\{v_{C},$ $C\in\mathbf{C}(p)\}$ form a basis of $V_{\infty
}(1^{p}).$
\end{lemma}

\begin{proof}
By Lemma \ref{lem_relation}, each vector of $V_{\infty}(1^{p})$ can be
decomposed into a linear combination of vectors $v_{C}.$ Given a finite
family $\{v_{i}\}_{i\in I}$ of vectors $v_{C},$ there exists an integer $n$
such that all the vectors in $\{v_{i}\}_{i\in I}$ belongs to $%
V_{n}(1^{p}).\; $By Lemma 3.1.2 of \cite{lec3}, the vectors of $%
\{v_{i}\}_{i\in I}$ are linearly independent in $V_{n}(1^{p}),$ thus also in 
$V_{\infty}(1^{p})$.
\end{proof}

\bigskip

\noindent Consider $\lambda\in\mathcal{P}_{m}$ for any $p\in\{1,...,m\}$
denote by $\widehat{\lambda}_{p}$ the number of columns of height $p$ in the
Young diagram of $\lambda$.\ Set 
\begin{equation}
W_{\infty}(\lambda)=V_{\infty}(1)^{\otimes\widehat{\lambda}_{1}}\otimes
\cdot\cdot\cdot\otimes V_{\infty}(1^{m})^{\otimes\widehat{\lambda}_{m}}.
\label{W(lambda)}
\end{equation}
The natural basis of $W_{\infty}(\lambda)$ consists of the tensor products $%
v_{C_{r}}\otimes\cdot\cdot\cdot\otimes v_{C_{1}}$ of basis vectors $v_{C}$
of the previous section appearing in (\ref{W(lambda)}). For any positive
integer $n\geq m$ set $u_{p}^{(n)}=v_{\overline{n-p+1}}\wedge\cdot\cdot\cdot%
\wedge v_{\overline{n}}.\;$The vector $u_{p}^{(n)}$ is of highest weight $%
\phi _{n}(1^{p})$ in the restriction of $V_{\infty}(1^{p})$ to $U_{q}(\frak{g%
}_{n})$. In the restriction of $W_{\infty}(\lambda)$ to $U_{q}(\frak{g}%
_{n}), $ the vector $u_{\lambda}^{(n)}=(u_{1}^{(n)})^{\otimes\widehat{\lambda%
}_{1}}\otimes\cdot\cdot\cdot\otimes(u_{m}^{(n)})^{\otimes\widehat{\lambda}%
_{m}}$ is the highest weight vector of an irreducible component $%
V(u_{\lambda}^{(n)})$ isomorphic to $V_{n}(\lambda)$. By Lemma \ref
{lem-trans}, we have $V(u_{\lambda}^{(n)})\subset V(u_{\lambda}^{(n+1)})$.
By Remark $\mathrm{(i)}$ after Proposition \ref{prop_irre}, we obtain:

\begin{proposition}
The $U_{q}(\frak{g}_{\infty})$-modules $\cup_{n\geq m}V(u_{\lambda}^{(n)})$
and $V_{\infty}(\lambda)$ are isomorphic.
\end{proposition}

\section{Crystal basis and canonical basis for the modules in $\mathcal{M}$}

\subsection{Crystal basis}

Consider $M$ a module in the category $\mathcal{M}$. Since $M$ is
decomposable, there exists a finite number of irreducible $U_{q}(\frak{g}%
_{\infty})$-modules $V_{\infty}^{(1)},...,V_{\infty}^{(r)}$ in $\mathcal{M}$
such that 
\begin{equation*}
M=\bigoplus_{k=1}^{r}V_{\infty}^{(k)}.
\end{equation*}
Each irreducible module $V_{\infty}^{(k)}$ in $\mathcal{M}$ being isomorphic
to a module $V_{\infty}(\lambda)$ with $\lambda\in\mathcal{P}$, its
restriction to $U_{q}(\frak{g}_{n})$ contains a unique irreducible component 
$V_{n}^{(k)}$ isomorphic to $V_{n}(\lambda).$ Then we set $%
M_{n}=\oplus_{1\leq k\leq r}V_{n}^{(k)}$. We have $M_{n}\subset M_{n+1}$ and 
$M=\cup_{n\geq1}M_{n}.$

\bigskip

\noindent\textbf{Remark: }Consider $M$ and $N$ two $U_{q}(\frak{g}_{\infty})$%
-modules.\ For types $B,C$ and $D$ $\left( M\otimes N\right) _{n}$ is a
submodule of $M_{n}\otimes N_{n}$ but we have $\left( M\otimes N\right)
_{n}\neq M_{n}\otimes N_{n}$ in general. When $M=V_{\infty}(\lambda)$ and $%
N=V_{\infty}(\mu)$ are irreducible, $\left( V_{\infty}(\lambda)\otimes
V_{\infty}(\nu)\right) _{n}$ is the submodule of $V_{n}(\lambda)\otimes
V_{n}(\nu)$ obtained as the direct sum of its irreducible components of
highest weight $\phi_{n}(\nu)$ where $\left| \nu\right| =\left|
\lambda\right| +\left| \mu\right| $ (see Corollary \ref{dec_tens}).

\begin{definition}
\label{def_cb}Consider a $U_{q}(\frak{g}_{\infty})$-module $M\in\mathcal{M}$
and denote by $m$ the minimal positive integer such that $M_{m}\neq\{0\}.$ A
crystal basis for $M$ is a pair $(L_{\infty},B_{\infty})=$ $%
(L_{n},B_{n})_{n\geq m}$ verifying the following conditions:

\begin{enumerate}
\item  for any $n\geq1,$ $(L_{n},B_{n})$ is a crystal basis of $M_{n}.$

\item  for any $n\geq1,$ we have $L_{n}\subset L_{n+1}$ and $b+qL_{n}\in
B_{n}\Longrightarrow b+qL_{n+1}\in B_{n+1}.$

\item  $L_{\infty}=\cup_{n\geq m}L_{n}$ and $B_{\infty}=\{b+qL_{\infty}\in
L_{\infty}/qL_{\infty}\mid\exists n\geq m,b+qL_{n}\in B_{n}\}.$
\end{enumerate}
\end{definition}

\begin{lemma}
\label{lem_tec2}Consider a vector $v\in V_{n}(\lambda)$ of weight $\delta$
such that $\dim V_{n}(\lambda)_{\delta}=1$.\ Then there exists a unique
crystal basis $(L_{n}(\lambda),B_{n}(\lambda))$ of $V_{n}(\lambda)$ such
that $v\in L_{n}(\lambda)$ and $v+qL_{n}(\lambda)\in B_{n}(\lambda)$.
\end{lemma}

\begin{proof}
Let $u_{\lambda}$ be a highest weight vector of $V_{n}(\lambda)$ and $%
(L_{n}(u_{\lambda}),B_{n}(u_{\lambda}))$ the crystal basis corresponding to $%
u_{\lambda}$ by (\ref{bc}). Since $\dim V_{n}(\lambda)_{\delta}=1$, there
exists a unique vertex $\widetilde{v}+qL_{n}(u_{\lambda})$ of weight $\delta$
in the crystal $B_{n}(u_{\lambda})$.\ Moreover we must have $\widetilde
{v}=K(q)v$ with $K(q)\in\mathbb{Q}(q).$ Then $v_{\lambda}=\frac{1}{K(q)}%
u_{\lambda}$ is a highest weight vector in $V_{n}(\lambda)$ and the crystal
basis $(L_{n}(\lambda),B_{n}(\lambda))$ obtained from $v_{\lambda}$ as in (%
\ref{bc}) is such that $v\in L_{n}(\lambda)$ and $v+qL_{n}(\lambda)\in
B_{n}(\lambda).$ This shows the ``existence'' part of the lemma.

\noindent Now suppose that $v_{\lambda}$ and $v_{\lambda}^{\prime}$ are two
highest weight vectors in $V_{n}(\lambda)$ defining two crystal bases $%
(L_{n}(\lambda),B_{n}(\lambda))$ and $(L_{n}^{\prime}(\lambda),B_{n}^{\prime
}(\lambda))$ verifying the assertion of the lemma. Set $v_{\lambda}^{\prime
}=K(q)v_{\lambda}$ with $K(q)\in\mathbb{Q}(q)$. Let $\widetilde{f}$ $=%
\widetilde{f}_{i_{1}}\cdot\cdot\cdot\widetilde{f}_{i_{r}}$ be a path in $%
B_{n}(\lambda)$ joining the highest weight vertex $b_{\lambda}=v_{\lambda
}+qL_{n}(\lambda)$ to $v+qL_{n}(\lambda)$. Then $\widetilde{f}(v_{\lambda
})\neq0$ in $V_{n}(\lambda)$ and has weight $\delta$.\ Write $\widetilde
{f}(v_{\lambda})=F(q)v$ with $F(q)\in\mathbb{Q}(q)$.\ We have $\widetilde
{f}(v_{\lambda})\in L_{n}(\lambda),$ $v\in L_{n}(\lambda)$ and $v\notin
qL_{n}(\lambda)$ (otherwise $v+qL_{n}(\lambda)=0).$ This implies that $%
F(q)\in\mathbf{A}(q)$. Moreover, we must have $\widetilde{f}(v_{\lambda
})\equiv v\func{mod}qL_{n}(\lambda)$.\ Thus $F(0)=1.$ Similarly, we can
write $\widetilde{f}(v_{\lambda}^{\prime})=F^{\prime}(q)v$ with $%
F^{\prime}(q)\in\mathbf{A}(q)$ such that $F^{\prime}(0)=1$. By using the
equality $v_{\lambda}^{\prime}=K(q)v_{\lambda},$ we derive $\widetilde
{f}(v_{\lambda}^{\prime})=F^{\prime}(q)v=K(q)F(q)v$. Hence $K(q)=\frac
{F^{\prime}(q)}{F(q)}\in\mathbf{A}(q)$ and is such that $K(0)=1$. By Lemma 
\ref{lm_tec}, this implies that the crystal bases $(L_{n}(\lambda
),B_{n}(\lambda))$ and $(L_{n}^{\prime}(\lambda),B_{n}^{\prime}(\lambda))$
coincide.
\end{proof}

\begin{theorem}
\label{th_un_bc}Consider $\lambda\in\mathcal{P}_{m}$ and $v_{\lambda}^{(m)}$
a fixed vector of highest weight $\phi_{m}(\lambda)$ in $V_{m}(\lambda).$
Then $V_{\infty}(\lambda)$ contains a unique crystal basis $%
(L_{\infty}(\lambda),B_{\infty}(\lambda))=(L_{n}(\lambda),B_{n}(\lambda))_{n%
\geq m}$ such that for any $n\geq m$, $v_{\lambda}^{(m)}\in L_{n}(\lambda)$
and $v_{\lambda }^{(m)}+qL_{n}(\lambda)\in B_{n}(\lambda)$.
\end{theorem}

\begin{proof}
By the above lemma, for any integer $n\geq m,$ there exists a unique crystal
basis $(L_{n}(\lambda),B_{n}(\lambda))$ of $V_{n}(\lambda)$ such that $%
v_{\lambda}^{(m)}\in L_{n}(\lambda)$ and $v_{\lambda}^{(m)}+qL_{n}(\lambda)%
\in B_{n}(\lambda)$. We are going first to establish that $%
(L_{\infty}(\lambda),B_{\infty}(\lambda))=(L_{n}(\lambda),B_{n}(\lambda))_{n%
\geq m}$ is a crystal basis of $V_{\infty}(\lambda)$ in the sense of
Definition \ref{def_cb}.

\noindent Consider an integer $n\geq m.\;$Set $L_{n}^{\prime}(\lambda
)=L_{n+1}(\lambda)\cap V_{n}(\lambda)$ and 
\begin{equation*}
B_{n}^{\prime}(\lambda)=\{b+qL_{n}^{\prime}(\lambda)\mid
b+qL_{n+1}(\lambda)\in B_{n+1}(\lambda),b\in L_{n}^{\prime}(\lambda)\}.
\end{equation*}
Then $(L_{n}^{\prime}(\lambda),B_{n}^{\prime}(\lambda))$ is a crystal basis
of $V_{n}(\lambda)$ such that $v_{\lambda}^{(m)}\in L_{n}^{\prime}(\lambda)$
and $v_{\lambda}^{(m)}+qL_{n}^{\prime}(\lambda)\in B_{n}^{\prime}(\lambda)$.
By the previous lemma, we must have $(L_{n}^{\prime}(\lambda),B_{n}^{\prime
}(\lambda))=(L_{n}(\lambda),B_{n}(\lambda))$.\ Thus $L_{n}(\lambda)\subset
L_{n+1}(\lambda)$ and $b+qL_{n}(\lambda)\in B_{n}(\lambda)\Longrightarrow
b+qL_{n+1}(\lambda)\in B_{n+1}(\lambda).$ This shows that $%
(L_{\infty}(\lambda),B_{\infty}(\lambda))=(L_{n}(\lambda),B_{n}(\lambda))_{n%
\geq m}$ is a crystal basis of $V_{\infty}(\lambda)$.

\noindent Now if $(L_{\infty}(\lambda),B_{\infty}(\lambda))=(L_{n}(%
\lambda),B_{n}(\lambda))_{n\geq m}$ and $(L_{\infty}^{\prime}(\lambda
),B_{\infty}^{\prime}(\lambda))=(L_{n}^{\prime}(\lambda),B_{n}^{\prime
}(\lambda))_{n\geq m}$ are two crystal bases of $V_{\infty}(\lambda)$
verifying the assertion of the theorem, we deduce immediately from Lemma \ref
{lem_tec2} that for any $n\geq m,$ the crystal bases $(L_{n}(\lambda),B_{n}(%
\lambda))$ and $(L_{n}^{\prime}(\lambda),B_{n}^{\prime }(\lambda))$ coincide.
\end{proof}

\bigskip

\noindent\textbf{Remarks:}

\noindent$\mathrm{(i)}$\textrm{:} If $(L_{\infty},B_{\infty})$ is a crystal
basis of $M,$ $L_{\infty}$ is a crystal lattice for $M$ and $B$ a $\mathbb{Q}
$ basis of $L_{\infty}/qL_{\infty}$ compatible with the weight decomposition
of $M$ and stable under the action of the Kashiwara operators (see \cite{Kan}
Definition 4.2.3).\ Unfortunately, these conditions are not sufficient to
guarantee the unicity of the crystal basis of $V_{\infty}(\lambda)$ because
this module has no highest weight vector.

\noindent$\mathrm{(ii)}$\textrm{:} The crystal basis $(L_{\infty}(%
\lambda),B_{\infty}(\lambda))=(L_{n}(\lambda),B_{n}(\lambda))_{n\geq m}$ of
Theorem \ref{th_un_bc} verifies 
\begin{equation*}
L_{n}(\lambda)=L_{\infty}(\lambda)\cap V_{n}(\lambda)\text{ and }%
B_{n}(\lambda)=\{v+qL_{n}(\lambda)\mid v+qL_{\infty}(\lambda)\in
B_{\infty}(\lambda),v\in L_{n}(\lambda)\}.
\end{equation*}
For any positive integer $p$, the crystal basis of $V_{\infty}(1^{p})$ can
be explicitly described by using the $q$-wedge product realization given in 
\ref{subsec_qwed}. Consider the $\mathbf{A}(q)$-lattice $L_{\infty}^{(p)}=%
\oplus_{C\in\mathbf{C}(p)}\mathbf{A}(q)v_{C}$ in $V_{\infty}(1^{p})$ and set 
$B_{\infty}^{(p)}=\{v_{C}+qL^{(p)}\mid C\in\mathbf{C}(p)\}.$

\begin{proposition}
The crystal basis $(L_{\infty}(1^{p}),B_{\infty}(1^{p}))$ of $%
V_{\infty}(1^{p})$ associated to the highest weight vector $v^{(p)}=v_{%
\overline{p}}\wedge v_{\overline{p-1}}\wedge\cdot\cdot\cdot\wedge v_{%
\overline{1}}\in V_{p}(1^{p})$ by Theorem \ref{bc} verifies $%
L_{\infty}(1^{p})=L_{\infty}^{(p)}$ and $B_{\infty}(1^{p})=B_{\infty}^{(p)}.$
\end{proposition}

\begin{proof}
For any integer $n\geq1,$ one verifies (see Lemma 3.1.3 in \cite{lec3}) that 
$L_{n}((1^{p}))=L_{\infty}^{(p)}\cap V_{n}(1^{p})$ and $B_{n}(1^{p})=%
\{v+qL_{n}(1^{p})\mid v+qL^{(p)}\in B_{\infty}^{(p)}\}$ yields a crystal
basis of $V_{n}(1^{p})$ such that $v^{(p)}\in L_{n}(1^{p})$ and $%
v^{(p)}+qL_{n}(1^{p})\in B_{\infty}^{(p)}$. Moreover $L_{n}(1^{p})\subset
L_{n+1}(1^{p})$ and $b+qL_{n}(1^{p})\in B_{n}(\lambda)\Longrightarrow
b+qL_{n+1}(1^{p})\in B_{n+1}(1^{p}).$ Thus by unicity of the crystal basis
established in Theorem \ref{th_un_bc}, $L_{\infty}(1^{p})=\cup_{n\geq
p}L_{n}(1^{p})=L_{\infty}^{(p)}$ and $B_{\infty}(1^{p})=B_{\infty}^{(p)}$.
\end{proof}

\begin{proposition}
\label{prop_sdpt}Suppose that $(L,B)$ and $(L^{\prime},B^{\prime})$ are
crystal bases of the $U_{q}(\frak{g}_{\infty})$-modules $M$ and $M^{\prime}$
belonging to $\mathcal{M}$. Then

\begin{enumerate}
\item  $(L\oplus L^{\prime},B\cup B^{\prime})$ is a crystal basis of $%
M\oplus M^{\prime}$

\item  $(L\otimes L^{\prime},B\otimes B^{\prime})$ with $B\otimes B^{\prime
}=\{b\otimes b^{\prime}\mid b\in B,b^{\prime}\in B^{\prime}\}$ is a crystal
basis of $M\otimes M^{\prime}$ where the action of the Kashiwara operators $%
\widetilde{e}_{i}$ and $\widetilde{f}_{i}$ $i\in I$ on $B\otimes B^{\prime}$
is given by (\ref{TENS1}) and (\ref{TENS2}).
\end{enumerate}
\end{proposition}

\begin{proof}
The proof of $1$ is the same as in the finite rank case (see Theorem 4.2.10
in \cite{Kan}).

\noindent For $2,$ consider the sequence $(L_{n}^{\otimes},B_{n}^{\otimes
})_{n\geq m}$ where $L_{n}^{\otimes}=(L\otimes L^{\prime})\cap(M\otimes
M^{\prime})_{n}$ and 
\begin{equation*}
B_{n}^{\otimes}=\{(v\otimes v^{\prime}+qL_{n}^{\otimes}\mid v\otimes
v^{\prime}\in L_{n}^{\otimes},(v+qL)\otimes(v+qL^{\prime})\in B\otimes
B^{\prime}\}.
\end{equation*}
Each pair $(L_{n}^{\otimes},B_{n}^{\otimes})$ is a crystal basis of $%
(M\otimes M^{\prime})_{n}$ \ and the sequence $(L_{n}^{\otimes},B_{n}^{%
\otimes})_{n\geq m}$ verifies the conditions of Definition \ref{def_cb}. For
any integer $n$, the action of the Kashiwara operators on $B_{n}^{\otimes}$
is given by (\ref{TENS1}) and (\ref{TENS2}). We have $L\otimes
L^{\prime}=\cup_{n\geq m}L_{n}^{\otimes}$ and 
\begin{equation*}
B\otimes B^{\prime}=\{v\otimes v^{\prime}+qL\otimes L^{\prime}\mid\exists
n\geq m,v\otimes v^{\prime}+qL_{n}^{\otimes}\in B_{n}^{\otimes}\}.
\end{equation*}
Moreover, the conditions $L_{n}^{\otimes}\subset L_{n+1}^{\otimes}$ and $%
v+qL_{n}^{\otimes}\in B_{n}\Longrightarrow v+qL_{n+1}^{\otimes}\in B_{n+1}$
imply that the action of the Kashiwara operators $\widetilde{e}_{i},%
\widetilde{f}_{i}$ with $i\in I_{n}$ on $(L_{n},B_{n})$ and $%
(L_{n+1},B_{n+1})$ coincide. Thus the action of these operators on $B\otimes
B^{\prime}$ is also given by rules (\ref{TENS1}) and (\ref{TENS2}).
\end{proof}

\bigskip

\noindent Consider a $U_{q}(\frak{g}_{\infty})$-modules $M$ belonging to $%
\mathcal{M}$ with crystal basis $(L,B)$.\ Write $M=\oplus V_{\infty}^{(k)}$
for its decomposition into irreducible components.\ For any $k,$ let $%
(L^{(k)},B^{(k)})$ be the crystal basis of the irreducible module $V_{\infty
}^{(k)}.\;$By Proposition \ref{prop_sdpt}, we know that $(\oplus
L^{(k)},\cup B^{(k)})$ is also a crystal basis of $M$.

\begin{theorem}
\label{th_iso_bc}There exists an isomorphism of crystal bases 
\begin{equation*}
\psi:(L,B)\overset{\simeq}{\rightarrow}(\oplus L^{(k)},\cup B^{(k)}).
\end{equation*}
\end{theorem}

\begin{proof}
We write for short $V_{n}^{(k)}$ for $(V_{\infty}^{(k)})_{n}$. We have $%
M=\oplus V_{\infty}^{(k)},$ thus for any integer $n\geq1$, $M_{n}=\oplus
V_{n}^{(k)}.$ Set $L_{n}=L\cap M_{n}$ and $B_{n}=\{v+qL_{n}\mid v\in
L_{n},v+qL\in B\}$. Then $(L_{n},B_{n})$ is a crystal basis of $M_{n}.$
Similarly for any $k,$ $L_{n}^{(k)}=L^{(k)}\cap V_{n}^{(k)}$ and $%
B_{n}^{(k)}=\{v+qL_{n}^{(k)}\mid v\in L_{n}^{(k)},v+qL^{(k)}\in B^{(k)}\}$
define a crystal basis of $V_{n}^{(k)}.$ Thus $(\oplus L_{n}^{(k)},\cup
B_{n}^{(k)})$ is also a crystal basis of $M_{n}$.\ This implies that for any
integer $n\geq1,$ there exists an isomorphism of crystal bases 
\begin{equation*}
\psi_{n}:(L_{n},B_{n})\overset{\simeq}{\rightarrow}(\oplus L_{n}^{(k)},\cup
B_{n}^{(k)}).
\end{equation*}
The conditions 
\begin{equation*}
\left\{ 
\begin{array}{l}
L_{n}=L_{n+1}\cap M_{n},B_{n}=B_{n+1}\cap(L_{n}/qL_{n}) \\ 
L_{n}^{(k)}=L_{n+1}^{(k)}\cap V_{n}^{(k)},B_{n}^{(k)}=B_{n+1}^{(k)}\cap
(L_{n}^{(k)}/qL_{n}^{(k)})
\end{array}
\right. \text{ for any }n\geq1
\end{equation*}
shows that the restriction of $\psi_{n+1}$ to $(L_{n},B_{n})$ is also an
isomorphism of crystal bases from $(L_{n},B_{n})$ to $(\oplus
L_{n}^{(k)},\cup B_{n}^{(k)})$.\ Indeed, $\psi_{n+1}(L_{n})=\oplus
L_{n+1}^{(k)}\cap M_{n}=\oplus L_{n}^{(k)}$ for $L_{n}=L_{n+1}\cap M_{n}$
and $L_{n}^{(k)}=L_{n+1}^{(k)}\cap V_{n}^{(k)}$. Similarly 
\begin{equation*}
\psi_{n+1}(B_{n})=(\cup B_{n+1}^{(k)})\cap(\oplus L_{n}^{(k)}/\oplus
qL_{n}^{(k)})=\cup(B_{n+1}^{(k)}\cap L_{n}^{(k)}/qL_{n}^{(k)})=\cup
B_{n}^{(k)}.
\end{equation*}
This permits to choose the isomorphisms $\psi_{n}$ so that for any integer $%
n,$ $\psi_{n}$ is the restriction of $\psi_{n+1}$ to $(L_{n},B_{n})$. Then
one defines 
\begin{equation*}
\psi:(L,B)\rightarrow(\oplus L^{(k)},\cup B^{(k)})
\end{equation*}
by requiring that $\psi(v)=\psi_{n}(v)$ for any $v\in L_{n}$. The map $%
\psi:L\rightarrow\oplus L^{(k)}$ is a $\mathbf{A}(q)$-linear isomorphism
which commutes with the Kashiwara operators. Moreover, the induced $\mathbb{Q%
}$-linear isomorphism from $L/qL$ to $\oplus L^{(k)}/q\oplus L^{(k)}$
defines an isomorphism of crystals from $B$ to $\cup B^{(k)}$.\ Thus $\psi$
is an isomorphism of crystal bases.
\end{proof}

\bigskip

\noindent\textbf{Remark: }We deduce from the previous theorem that the
decomposition of $M$ into its irreducible components is given by the
decomposition of $B$ into its connected components.

\subsection{Canonical basis for the modules $V_{\infty}(\protect\lambda)$}

Denote by $F\mapsto\overline{F}$ the involution of $U_{q}(\frak{g}_{\infty})$
defined as the ring automorphism satisfying 
\begin{equation*}
\overline{q}=q^{-1},\text{ \ }t_{i}=t_{i}^{-1},\text{ \ \ }\overline{e_{i}}%
=e_{i},\text{ \ \ }\overline{f_{i}}=f_{i}\text{ \ \ \ for }i\in I.
\end{equation*}
Let $U_{\mathbb{Q}}^{-}$ (resp.\ $U_{n,\mathbb{Q}}^{-}$) be the subalgebra
of $U_{q}(\frak{g}_{\infty})$ generated over $\mathbb{Q}[q,q^{-1}]$ by the $%
f_{i}^{(k)},i\in I$ (resp. $i\in I_{n}$).\ Let us first recall the
definition of the canonical basis of $V_{n}(\lambda).$ Let $%
v_{\lambda}^{n}\in V_{n}(\lambda)$ be a highest weight vector of weight $%
\phi_{n}(\lambda)$.\ By writing each vector $v$ of $V_{n}(\lambda)$ in the
form $v=Fv_{\lambda}^{(n)}$ where $F\in U_{q}^{-}(\frak{g}_{n})$, we obtain
an involution $I_{n}$ of $V_{n}(\lambda)$ defined by 
\begin{equation}
I_{n}(v)=\overline{F}v_{\lambda}^{(n)}.  \label{defI}
\end{equation}
Set $V_{n,\mathbb{Q}}(\lambda)=U_{n,\mathbb{Q}}^{-}v_{\lambda}^{(n)}$ and
denote by $(L_{n}(\lambda),B_{n}(\lambda))$ the crystal basis of $%
V_{n}(\lambda)$ determined by $v_{\lambda}^{(n)}.$

\begin{theorem}
\label{TH_K_2}(Kashiwara) There exists a unique $\mathbb{Q[}q,q^{-1}]$-basis 
$G_{n}(\lambda)=\{G_{n}(b)\mid b\in B_{n}(\lambda)\}$ of $V_{n,\mathbb{Q}%
}(\lambda)$ such that: 
\begin{gather}
G_{n}(b)\equiv b\text{ }\mathrm{mod}\text{ }qL_{n}(\lambda),
\label{cond_cong} \\
I_{n}(G_{n}(b))=G_{n}(b).  \label{cond_invo}
\end{gather}
\end{theorem}

\noindent The basis $G_{n}(\lambda)$ is called the canonical basis of $%
V_{n}(\lambda)$. It is determined uniquely up to the choice of a highest
weight vector in $V_{n}(\lambda)$.

\noindent\textbf{Remark: }Suppose that the highest weight vector $v_{\lambda
}^{(n)}$ is replaced by $(v_{\lambda}^{(n)})^{\prime}=K(q)v_{\lambda}^{(n)}$
where $K(q)\in\mathbb{Q}(q)$. Then the canonical basis $G_{n}^{^{\prime}}(%
\lambda)$ defined from $(v_{\lambda}^{(n)})^{\prime}$ verifies $%
G_{n}^{\prime}(\lambda)=K(q)G_{n}(b)$. Moreover, we have $I^{\prime}(v)=%
\frac{K(q)}{K(q^{-1})}I(v)$ where $I_{n}^{\prime}$ is the involution
corresponding to $(v_{\lambda}^{(n)})^{\prime}$.

\begin{lemma}
\label{lem_tech3}Consider a vector $v\in V_{n}(\lambda)$ of weight $\delta$
such that $\dim V_{n}(\lambda)_{\delta}=1$.\ Then there exists a unique
highest weight vector $v_{\lambda}^{(n)}\in V_{n}(\lambda)$ such that $v$ is
a vector of the canonical basis $G_{n}(\lambda)$ obtained from $v_{\lambda
}^{(n)}$.
\end{lemma}

\begin{proof}
Let $u_{\lambda}$ be a highest weight vector of $V_{n}(\lambda)$. Write $%
G_{n}^{u_{\lambda}}(\delta)$ for the unique vector of weight $\delta$ in the
canonical basis determined by $u_{\lambda}$. We must have $%
G_{n}^{u_{\lambda}}(\delta)=K(q)v$ where $K(q)\in\mathbb{Q}(q)$ for $\dim
V_{n}(\lambda )_{\delta}=1$.\ Set $v_{\lambda}^{(n)}=\frac{1}{K(q)}%
u_{\lambda}.$ The vector of weight $\delta$ in the canonical basis
determined by $v_{\lambda}^{(n)}$ is equal to $\frac{1}{K(q)}%
G_{n}^{u_{\lambda}}(\delta)=v$ (see remark above) as required. This shows
the ``existence'' part of the lemma.

\noindent Now suppose that $v_{\lambda}^{(n)}$ and $(v_{\lambda}^{(n)})^{^{%
\prime}}$ are two highest weight vectors in $V_{n}(\lambda)$ such that $v$
belongs to the canonical bases $G_{n}(\lambda)$ and $G_{n}^{\prime
}(\lambda) $ determined respectively by $v_{\lambda}^{(n)}$ and $(v_{\lambda
}^{(n)})^{^{\prime}}$. Write $(v_{\lambda}^{(n)})^{\prime}=K(q)v_{\lambda
}^{(n)}$ where $K(q)\in\mathbb{Q}(q)$. Denote by $G_{n}(\delta)$ and $%
G_{n}^{\prime}(\delta)$ the vectors of weight $\delta$ respectively in the
canonical bases $G_{n}(\lambda)$ and $G_{n}^{\prime}(\lambda).$ We must have 
$G_{n}^{\prime}(\delta)=K(q)G_{n}(\delta).$ Since $\dim V_{n}(\lambda
)_{\delta}=1$ and $v$ belongs to $G_{n}(\lambda)\cap
G_{n}^{\prime}(\lambda), $ we obtain $G_{n}(\delta)=G_{n}^{\prime}(\delta)=v.
$ Thus $K(q)=1$ and $v_{\lambda}^{(n)}=(v_{\lambda}^{(n)})^{^{\prime}}$.
\end{proof}

\bigskip

\noindent The canonical basis $G_{n+1}(\lambda)$ defined by $v_{\lambda
}^{(n+1)}$ contains a unique vector $u_{\lambda}^{(n)}$ of weight $\phi
_{n}(\lambda)$. This vector is of highest weight in $V_{n}(\lambda).$ Denote
by $G_{n}^{(n+1)}(\lambda)$ the canonical basis of $V_{n}(\lambda)$ defined
from $u_{\lambda}^{(n)}$. Write $I_{n}^{(n+1)}$ for the corresponding
involution defined on $V_{n}(\lambda)$ and set $V_{n,\mathbb{Q}%
}^{(n+1)}(\lambda)=U_{n,\mathbb{Q}}^{-}u_{\lambda}^{(n)}$.

\begin{lemma}
\label{lem_tech4}For any $v\in V_{n}(\lambda),$ we have $%
I_{n}^{(n+1)}(v)=I_{n+1}(v).$ Moreover $V_{n,\mathbb{Q}}^{(n+1)}(\lambda)%
\subset V_{n+1,\mathbb{Q}}(\lambda)$ and $G_{n}^{(n+1)}(\lambda)\subset
G_{n+1}(\lambda).$
\end{lemma}

\begin{proof}
The vector $u_{\lambda}^{(n)}$ belongs to $V_{n+1}(\lambda)$, thus there
exists $R\in U_{q}(\frak{g}_{n+1})$ such that $u_{\lambda}^{(n)}=R\cdot
v_{\lambda}^{(n+1)}$. Consider $v\in V_{n}(\lambda)$ and set $v=F\cdot
v_{\lambda}^{n}$ where $F\in U_{q}(\frak{g}_{n})$. Then $v=FR\cdot
v_{\lambda }^{(n+1)}.$ We have $I_{n+1}(v)=\overline{FR}\cdot
v_{\lambda}^{(n+1)}=\overline{F}\ (\overline{R}\cdot v_{\lambda}^{(n+1)})$.
Since $u_{\lambda }^{(n)}\in G_{n+1}(\lambda),$ we must have $%
I_{n+1}(u_{\lambda}^{(n)})=u_{\lambda}^{(n)}=\overline{R}\cdot
v_{\lambda}^{(n+1)}$. This gives $I_{n+1}(v)=\overline{F}\cdot
u_{\lambda}^{(n)}=I_{n}^{(n+1)}(v)$. The inclusion $V_{n,\mathbb{Q}%
}^{(n+1)}(\lambda)\subset V_{n+1,\mathbb{Q}}(\lambda)$ is immediate. We
deduce that the vectors $G_{n}^{(n+1)}(b)$ of $G_{n}^{(n+1)}(\lambda)$
belongs to $V_{n+1,\mathbb{Q}}(\lambda)$. Moreover they verify $%
G_{n}^{(n+1)}(b)\equiv b$ $\mathrm{mod}$ $qL_{n+1}(\lambda)$ and $%
I_{n+1}(G_{n}^{(n+1)}(b))=G_{n}^{(n+1)}(b).$ Thus $G_{n}^{(n+1)}(\lambda)%
\subset G_{n+1}(\lambda).$
\end{proof}

\bigskip

\noindent Now suppose that $\lambda\in\mathcal{P}_{m}$ and fix a vector $%
v_{\lambda}^{(m)}$ in $V_{m}(\lambda)$ of highest weight $\phi_{m}(\lambda)$%
. For any $n\geq m,$ $v_{\lambda}^{(m)}$ belongs to $V_{n}(\lambda)$ and $%
\dim V_{n}(\lambda)_{\phi_{m}(\lambda)}=1$. Thus by the above lemma, there
is a unique highest weight vector $v_{\lambda}^{(n)}\in V_{n}(\lambda)$ such
that $v_{\lambda}^{(m)}$ belongs to the canonical basis $G_{n}(\lambda)$
determined by $v_{\lambda}^{(n)}.$ This yields an involution $I_{n}$ on $%
V_{n}(\lambda)$ and a subspace $V_{n,\mathbb{Q}}(\lambda)$ which depends
only on the choice of $v_{\lambda}^{(m)}$.

\noindent For any integer $n,$ the family of canonical bases $%
(G_{n}(\lambda))_{n\geq1}$ determined by the condition $v_{\lambda}^{(m)}\in
G_{n}(\lambda)$ verify $G_{n}(\lambda)\subset G_{n+1}(\lambda).$ Indeed, $%
G_{n+1}(\lambda)$ contains a unique vector $u_{\lambda}^{(n)}$ of weight $%
\phi_{n}(\lambda).$ This vector is of highest weight in $V_{n}(\lambda)$. By
Lemma \ref{lem_tech4}, we obtain that $v_{\lambda}^{m}$ belongs to the
canonical bases of $V_{n}(\lambda)$ determined by $u_{\lambda}^{(n)}$. Thus
by Lemma \ref{lem_tech3} we have $u_{\lambda}^{(n)}=v_{\lambda}^{(n)}$. In
particular for any $n\geq m,$ $I_{n}$ is the restriction of $I_{n+1}$ on $%
V_{n}(\lambda).$ This permits to define the bar involution on $V_{\infty
}(\lambda)$ by 
\begin{equation*}
\overline{v}=I_{n}(v)\text{ for any }v\in V_{n}(\lambda)\text{.}
\end{equation*}
Set $V_{\mathbb{Q}}(\lambda)=\cup_{n\geq m}V_{n,\mathbb{Q}}(\lambda)$ and
consider the crystal basis $(L_{\infty}(\lambda),B_{\infty}(\lambda))$ of $%
V_{\infty}(\lambda)$ determined by $v_{\lambda}^{(m)}$ as in Theorem \ref
{th_un_bc}. The following theorem follows immediately from the previous
arguments.

\begin{theorem}
There exists a unique $\mathbb{Q[}q,q^{-1}]$-basis $G_{\infty}(\lambda
)=\{G_{\infty}(b)\mid b\in B_{\infty}(T)\}$ of $V_{\mathbb{Q}}(\lambda)$
such that: 
\begin{gather*}
G_{\infty}(b)\equiv b\text{ }\mathrm{mod}\text{ }qL_{\infty}(\lambda), \\
\overline{G_{\infty}(b)}=G_{\infty}(b).
\end{gather*}
\end{theorem}

\noindent We will call $G_{\infty}(\lambda)$ the canonical basis of $%
V_{\infty}(\lambda)$. It is determined uniquely up to the choice of a
highest weight vector in $V_{m}(\lambda)$.

\bigskip

\noindent\textbf{Remark: }One has $G_{\infty}(\lambda)=\cup_{n\geq
m}G_{n}(\lambda)$ where $G_{n}(\lambda)$ is the canonical basis of $%
V_{n}(\lambda)$ determined by $v_{\lambda}^{(m)}$.

\section{Robinson-Schensted correspondences and bi-crystals\label{sec_plac}}

\subsection{Insertion schemes\label{sub_sec_tab}}

\noindent We have seen in \ref{subsec_RS} that the vertices of $B_{n}(1)$
are labelled by the letters of the alphabet $\mathcal{X}_{n}$.\ This permits
to identify the vertices of the crystal graph $G_{n}=\underset{l\geq0}{%
\tbigoplus }B_{n}(1)^{\otimes l}$ with the words on $\mathcal{X}_{n}.\;$For
any $w\in G_{n}$ we have $\mathrm{wt}(w)=(d_{1},...,d_{n})\in P_{n}$ where
for all $i=1,...,n$ $d_{i}$ is the number of letters $\overline{i}$ of $w$
minus its number of letters $i.$

\noindent Consider $\lambda\in\mathcal{P}_{m}$ with $m\leq n$ and denote by $%
Y(\lambda)$ the Young diagram of shape $\lambda.\;$Write $T_{n,\lambda}$ for
the filling of $Y(\lambda)$ whose $k$-th row contains only letters $%
\overline{n-k+1}.$ Let $b_{n,\lambda}$ be the vertex of $B_{n}(1)^{\otimes
\left| \lambda\right| }$ obtained by column reading $T_{n,\lambda}$, that is
by reading $T_{n,\lambda}$ from right to left and top to bottom. Kashiwara
and Nakashima realize $B_{n}(\lambda)$ as the connected component of the
tensor power $B_{n}(1)^{\otimes\left| \lambda\right| }$ of highest weight
vertex $b_{n,\lambda}$. The Kashiwara-Nakashima tableaux of type $%
A_{n-1},B_{n},C_{n}$, $D_{n}$ and shape $\lambda$ are defined as the
fillings of $Y(\lambda)$ having a column reading which belongs to $%
B_{n}(\lambda)$ (recall that $B_{n}(\lambda)$ depends on the root system
considered). In the sequel we will denote by $\mathrm{w}(T)$ the column
reading of the tableau $T.$

\noindent Write $\mathbf{T}_{n}(\lambda)$ respectively for the sets of
Kashiwara-Nakashima's tableaux of shape $\lambda$. Set $\mathbf{T}_{n}=%
\underset{\lambda\in\mathcal{P}_{m},m\leq n}{\cup}\mathbf{T}_{n}(\lambda)$
and $\mathbf{T}=\cup_{n\geq1}\mathbf{T}_{n}$.\ In the sequel we only
summarize the combinatorial description of $\mathbf{T}^{n}(\lambda)$ and
refer the reader to \cite{KN}, \cite{Lec} and \cite{lec2} for more details.

\noindent When $\lambda$ is a column partition, the tableaux of $\mathbf{T}%
^{n}(\lambda)$ are called $n$-admissible columns. The $n$-admissible columns
of types $A_{n-1},B_{n},C_{n}$ and $D_{n}$ are in particular columns of
types $A_{n-1},B_{n},C_{n}$ and $D_{n}$ that is, columns of types $A_{\infty
},B_{\infty},C_{\infty}$ and $D_{\infty}$ (see (\ref{col})) with letters in $%
\mathcal{X}_{n}$.\ For type $A_{n-1},$ each column is $n$-admissible.\ This
is not true for the types $B_{n},C_{n}$ and $D_{n}.\;$More precisely a
column $C$ of (\ref{col}) is $n$-admissible if and only if it can be
duplicated following a simple algorithm described in \cite{lec2} into a pair 
$(lC,rC)$ of columns without pair of opposite letters $(x,\overline{x})$
(the letter $0$ is counted as the pair $(0,\overline{0})$) and containing
only letters $a$ such that $a\in\mathcal{X}_{n}.$

\begin{example}
For the column $C= 
\begin{tabular}{|l|}
\hline
$\mathtt{\bar{3}}$ \\ \hline
$\mathtt{\bar{1}}$ \\ \hline
$\mathtt{0}$ \\ \hline
$\mathtt{1}$ \\ \hline
$\mathtt{2}$ \\ \hline
\end{tabular}
$ of type $B$ we have $lC=$%
\begin{tabular}{|l|}
\hline
$\mathtt{\bar{5}}$ \\ \hline
$\mathtt{\bar{4}}$ \\ \hline
$\mathtt{\bar{3}}$ \\ \hline
$\mathtt{1}$ \\ \hline
$\mathtt{2}$ \\ \hline
\end{tabular}
and $rC=$%
\begin{tabular}{|l|}
\hline
$\mathtt{\bar{3}}$ \\ \hline
$\mathtt{\bar{1}}$ \\ \hline
$\mathtt{3}$ \\ \hline
$\mathtt{4}$ \\ \hline
$\mathtt{5}$ \\ \hline
\end{tabular}
. Hence $C$ is $5$-admissible but not $n$-admissible for $n\leq4.$
\end{example}

\noindent For a general $\lambda$ a tableau $T\in\mathbf{T}^{n}(\lambda)$
can be regarded as a filling of the Young diagram of shape $\lambda$ such
that

\begin{itemize}
\item  $T=C_{1}\cdot\cdot\cdot C_{r}$ where the columns $C_{i}$ of $T$ are $%
n $-admissible,

\item  for any $i\in\{1,...r-1\}$ the columns of the tableau $%
r(C_{i})l(C_{i+1})$ weakly increase from left to right and do not contain
special configurations (detailed in \cite{KN} and \cite{lec2}) when $T$ is
of type $D_{n}.$
\end{itemize}

\noindent\textbf{Remarks:}

\noindent$\mathrm{(i):}$\textbf{\ }We have $\mathbf{T}^{n}(\lambda )\subset%
\mathbf{T}^{n+1}(\lambda)$ since the $n$-admissible columns are also $(n+1)$%
-admissible and the duplication process of a column does not depend on $n.$

\noindent$\mathrm{(ii):}$ For type $A_{n-1},$ we have $lC=rC=C$ and $%
T^{n}(\lambda)$ is the set of semistandard tableaux of shape $\lambda$.

\bigskip

There exist insertion schemes related to each classical root system \cite
{Lec}, \cite{lec2} analogous for Kashiwara-Nakashima's tableaux to the well
known bumping algorithm on semistandard tableaux.

\noindent Denote by $\sim_{n}$ the equivalence relation defined on the
vertices of $G_{n}$ by $w_{1}\sim_{n}w_{2}$ if and only if $w_{1}$ and $%
w_{2} $ belong to the same connected component of $G_{n}.\;$For any word $w,$
the insertion schemes permit to compute the unique tableau $P_{n}(w)$ such
that $w\sim_{n}\mathrm{w}(P_{n}(w))$. In fact $\sim_{n}$ is a congruence $%
\equiv_{n}$ \cite{Lec} \cite{lec2} \cite{lit} which depend on the type
considered. For type $A_{n-1},$ $\equiv_{n}$ is the congruence defined by
the Knuth relations.\ For types $B_{n},C_{n}$ and $D_{n},$ it is obtained as
the quotient of the free monoid of words on $\mathcal{X}_{n}$ by two kinds
of relations.

\noindent The first consists of relations of length $3$ analogous to the
Knuth relations.\ In fact these relations are precisely those which are
needed to describe the insertion $x\rightarrow C$ of a letter $x$ in a $n$%
-admissible column $C= 
\begin{tabular}{|l|}
\hline
$a$ \\ \hline
$b$ \\ \hline
\end{tabular}
$ such that 
\begin{tabular}{|l|}
\hline
$a$ \\ \hline
$b$ \\ \hline
$x$ \\ \hline
\end{tabular}
is not a column. This can be written 
\begin{equation}
x\rightarrow 
\begin{tabular}{|l|}
\hline
$a$ \\ \hline
$b$ \\ \hline
\end{tabular}
= 
\begin{tabular}{c|c|}
\cline{2-2}
& $a$ \\ \hline
\multicolumn{1}{|c|}{$x$} & $b$ \\ \hline
\end{tabular}
= 
\begin{tabular}{|l|l}
\hline
$a^{\prime}$ & \multicolumn{1}{|l|}{$x^{\prime}$} \\ \hline
$b^{\prime}$ &  \\ \cline{1-1}
\end{tabular}
\label{trans_ele}
\end{equation}
and contrary to the insertion scheme for the semistandard tableaux the sets $%
\{a^{\prime},b^{\prime},x^{\prime}\}$ and $\{a,b,c\}$ are not necessarily
equal (i.e. the relations are not homogeneous in general).

\noindent Next we have the contraction relations which do not preserve the
length of the words.\ These relations are precisely those which are needed
to describe the insertion $x\rightarrow C$ of a letter $x$ such that $%
\overline{n}\leq x\leq n$ in a $n$-admissible column $C$ such that 
\begin{tabular}{|l|}
\hline
$C$ \\ \hline
$x$ \\ \hline
\end{tabular}
(obtained by adding the letter $x$ on bottom of $C$) is a column which is
not $n$-admissible. In this case 
\begin{tabular}{|l|}
\hline
$C$ \\ \hline
$x$ \\ \hline
\end{tabular}
is necessarily $(n+1)$-admissible and have to be contracted to give a $n$%
-admissible column. We obtain $x\rightarrow C=\widetilde{C}$ with $%
\widetilde{C}$ a $n$-admissible column of height $h(C)$ or $h(C)-1.$

\noindent The insertion of the letter $x$ in a $n$-admissible column $C$ of
arbitrary height such that 
\begin{tabular}{|l|}
\hline
$C$ \\ \hline
$x$ \\ \hline
\end{tabular}
is not a column can then be pictured by 
\begin{equation*}
x\rightarrow 
\begin{tabular}{|c|}
\hline
$a_{1}$ \\ \hline
$\cdot$ \\ \hline
$a_{k-2}$ \\ \hline
$a_{k-1}$ \\ \hline
$a_{k}$ \\ \hline
\end{tabular}
= 
\begin{tabular}{c|c|}
\cline{2-2}
\ \ \ \ \  & $a_{1}$ \\ \cline{2-2}
& $\cdot$ \\ \cline{2-2}
& $a_{k-2}$ \\ \cline{2-2}
& $a_{k-1}$ \\ \hline
\multicolumn{1}{|c|}{$x$} & $a_{k}$ \\ \hline
\end{tabular}
= 
\begin{tabular}{c|c}
\cline{2-2}
& \multicolumn{1}{|c|}{$a_{1}$} \\ \cline{2-2}
& \multicolumn{1}{|c|}{$\cdot$} \\ \cline{2-2}
& \multicolumn{1}{|c|}{$a_{k-2}$} \\ \hline
\multicolumn{1}{|c|}{$\delta_{k-1}$} & \multicolumn{1}{|c|}{$y$} \\ \hline
\multicolumn{1}{|c|}{$d_{k}$} &  \\ \cline{1-1}
\end{tabular}
=\cdot\cdot\cdot= 
\begin{tabular}{|c|c}
\hline
$d_{1}$ & \multicolumn{1}{|c|}{$\ \ z$ \ \ } \\ \hline
$\cdot$ &  \\ \cline{1-1}
$\cdot$ &  \\ \cline{1-1}
$d_{k-1}$ &  \\ \cline{1-1}
$d_{k}$ &  \\ \cline{1-1}
\end{tabular}
\end{equation*}
that is, one elementary transformation (\ref{trans_ele}) is applied to each
step. One proves that $x\rightarrow C$ is then a tableau of $\mathbf{T}%
^{(n)} $ with two columns respectively of height $h(C)$ and $1$.

\bigskip

\noindent Now we can define the insertion $x\rightarrow T$ of the letter $x$
such that $\overline{n}\leq x\leq n$ in the tableau $T\in\mathbf{T}%
^{n}(\lambda)$.\ Set $T=C_{1}\cdot\cdot\cdot C_{r}$ where $C_{i},$ $%
i=1,...,r $ are the $n$-admissible columns of $T.\;$

\begin{enumerate}
\item  When 
\begin{tabular}{|l|}
\hline
$C_{1}$ \\ \hline
$x$ \\ \hline
\end{tabular}
is not a column, write $x\rightarrow C= 
\begin{tabular}{|l|l|}
\hline
$C_{1}^{\prime}$ & $y$ \\ \hline
\end{tabular}
$ where $C_{1}^{\prime}$ is an admissible column of height $h(C_{1})$ and $y$
a letter.\ Then $x\rightarrow T=C_{1}^{\prime}(y\rightarrow C_{2}\cdot
\cdot\cdot C_{r})$ that is, $x\rightarrow T$ is the juxtaposition of $%
C_{1}^{\prime}$ with the tableau $\widehat{T}$ obtained by inserting $y$ in
the tableau $C_{2}\cdot\cdot\cdot C_{r}.$

\item  When 
\begin{tabular}{|l|}
\hline
$C_{1}$ \\ \hline
$x$ \\ \hline
\end{tabular}
is a $n$-admissible column, $x\rightarrow T$ is the tableau obtained by
adding a box containing $x$ on bottom of $C_{1}$.\ 

\item  When 
\begin{tabular}{|l|}
\hline
$C_{1}$ \\ \hline
$x$ \\ \hline
\end{tabular}
is a column which is not $n$-admissible, write $x\rightarrow C=\widetilde{C}$
and set $\mathrm{w}(\widetilde{C})=y_{1}\cdot\cdot\cdot y_{s}$ where the $%
y_{i}$'s are letters. Then $x\rightarrow
T=y_{s}\rightarrow(y_{s-1}\rightarrow(\cdot\cdot\cdot y_{1}\rightarrow%
\widehat{T}))$ that is $x\rightarrow T$ is obtained by inserting
successively the letters of $\widetilde{C}$ into the tableau $\widehat{T}%
=C_{2}\cdot\cdot\cdot C_{r}$. Note that there is no new contraction during
these $s$ insertions.
\end{enumerate}

\noindent\textbf{Remarks:}

\noindent$\mathrm{(i)}\mathbf{:}$ The $P_{n}$-symbol defined above can be
computed recursively by setting $P_{n}(w)= 
\begin{tabular}{|l|}
\hline
$w$ \\ \hline
\end{tabular}
$ if $w$ is a letter and $P_{n}(w)=x\rightarrow P_{n}(u)$ where $w=ux$ with $%
u$ a word and $x$ a letter otherwise.

\noindent$\mathrm{(ii)}\mathbf{:}$ Consider $T\in\mathbf{T}^{n}(\lambda
)\subset\mathbf{T}^{n+1}(\lambda)$ and a letter $x$ such that $\overline
{n}\leq x\leq n.$ The tableau obtained by inserting $x$ in $T$ may depend
wether $T$ is regarded as a tableau of $\mathbf{T}^{n}(\lambda)$ or as a
tableau of $\mathbf{T}^{n+1}(\lambda)$.\ Indeed if 
\begin{tabular}{|l|}
\hline
$C_{1}$ \\ \hline
$x$ \\ \hline
\end{tabular}
is not $n$-admissible then it is necessarily $(n+1)$-admissible since $C_{1}$
is $n$-admissible. Hence there is no contraction during the insertion $%
x\rightarrow T$ when it is regarded as a tableau of $\mathbf{T}%
^{n+1}(\lambda).$

\noindent$\mathrm{(iii)}\mathbf{:}$ Consider $w\in\mathcal{X}_{n}$.\ From $%
\mathrm{(ii)}\mathbf{\ }$we deduce that there exists an integer $N\geq n$
minimal such that $P_{N}(w)$ can be computed without using contraction
relation. Then for any $k\geq N,$ $P_{k}(w)=P_{N}(w).$

\noindent$\mathrm{(iv)}\mathbf{:}$ Similarly to the bumping algorithm for
semistandard tableaux, the insertion algorithms described above are
reversible.

\bigskip

\noindent We give below the plactic relations of length $3$. The contraction
relations will not be used in the sequel:

\subsubsection{For type $A_{n-1}$}

\begin{equation}
\left\{ 
\begin{tabular}{l}
$abx=bax$ if $a<x\leq b$ \\ 
$abx=axb$ if $x\leq a<b$%
\end{tabular}
\right.  \label{plA}
\end{equation}
where $a,b,x$ are letters of $\mathcal{X}_{A_{n-1}}.$

\subsubsection{For type $B_{n}$}

\begin{gather}
\left\{ 
\begin{tabular}{l}
$abx\equiv bax$ if $a<x<b$ and $b\neq\overline{a}$ \\ 
$abx\equiv axb$ if $x<a<b$ and $b\neq\overline{x}$%
\end{tabular}
\right. ,\left\{ 
\begin{tabular}{l}
$axx\equiv xax$ if $a<x$, $x\neq\overline{a}$ and $x\neq0$ \\ 
$xbx\equiv xxb$ if $x<b$, $b\neq\overline{x}$ and $x\neq0$%
\end{tabular}
\right. ,  \label{plB} \\
\left\{ 
\begin{tabular}{l}
$\overline{b}bx\equiv(b+1)(\overline{b+1})x$ if $\overline{b}\leq x\leq b$
and $b\neq0$ \\ 
$ab\overline{b}\equiv a(\overline{b-1})(b-1)$ if $\overline{b}<a<b$ and $%
(a,b)\neq(0,1)$%
\end{tabular}
\right. ,\left\{ 
\begin{tabular}{l}
$00x\equiv0x0\text{ if }x\leq\overline{1}$ \\ 
$0b0\equiv b00.$ if $b\geq1$%
\end{tabular}
\right. ,01\overline{1}\equiv1\overline{1}0  \notag
\end{gather}
where $a,b,x$ are letters of $\mathcal{X}_{B_{n}}.$

\subsubsection{For type $C_{n}$}

\begin{equation}
\left\{ 
\begin{tabular}{l}
$abx\equiv bax$ if $a<x\leq b$ and $b\neq\overline{a}$ \\ 
$abx\equiv axb$ if $x\leq a<b$ and $b\neq\overline{x}$%
\end{tabular}
\right. \text{ and }\left\{ 
\begin{tabular}{l}
$\overline{b}bx\equiv(b+1)(\overline{b+1})x$ if $\overline{b}\leq x\leq b$
\\ 
$ab\overline{b}\equiv a(\overline{b-1})(b-1)$ if $\overline{b}<a<b$%
\end{tabular}
\right.  \label{plC}
\end{equation}
where $a,b,x$ are letters of $\mathcal{X}_{C_{n}}.$

\subsubsection{For type $D_{n}$}

\begin{gather}
\left\{ 
\begin{tabular}{l}
$abx\equiv bax$ if $a<x\leq b$ and $b\neq\overline{a}$ \\ 
$abx\equiv axb$ if $x\leq a<b$ and $b\neq\overline{x}$%
\end{tabular}
\right. ,\left\{ 
\begin{tabular}{l}
$\overline{b}bx\equiv(b+1)(\overline{b+1})x$ if $\overline{b}\leq x\leq b$
\\ 
$ab\overline{b}\equiv a(\overline{b-1})(b-1)$ if $\overline{b}<a<b$%
\end{tabular}
\right.  \label{plD} \\
\left\{ 
\begin{tabular}{l}
$1\,b\overline{1}\equiv b1\,\overline{1}$ if $b\geq2$ \\ 
$\overline{1}\,b\,1\equiv b\,\overline{1}1$ if $b\geq2$%
\end{tabular}
\right. ,\text{ }\left\{ 
\begin{tabular}{l}
$1\overline{1}\,x\equiv1x\,\overline{1}$ if $x\leq\overline{2}$ \\ 
$\overline{1}1x\equiv\overline{1}\,x1$ if $x\leq\overline{2}$%
\end{tabular}
\right. ,\left\{ 
\begin{tabular}{l}
$\overline{1}\,11\equiv2\overline{2}1$ \\ 
$1\overline{1}\,\overline{1}\equiv2\overline{2}\,\overline{1}$%
\end{tabular}
\right. ,\text{ }\left\{ 
\begin{tabular}{l}
$12\overline{1}\equiv11\overline{1}$ \\ 
$\overline{1}21\equiv\overline{1}\,\overline{1}1$%
\end{tabular}
\right.  \notag
\end{gather}
where $a,b,x$ are letters of $\mathcal{X}_{D_{n}}.$

\subsection{Plactic monoids for types $A_{\infty},B_{\infty},C_{\infty}$ and 
$D_{\infty}$}

Consider $\lambda\in\mathcal{P}_{m}$, $v_{\lambda}^{(m)}\in V_{m}(\lambda)$
a highest weight vector and $(L_{\infty}(\lambda),B_{\infty}(\lambda))$ the
crystal basis of $V_{\infty}(\lambda)$ determined by $v_{\lambda}^{(m)}$. By
definition of the crystal basis $(L_{\infty}(\lambda),B_{\infty}(\lambda)),$
there is a natural embedding of crystals 
\begin{equation*}
\iota_{n}:\left\{ 
\begin{tabular}{l}
$B_{n}(\lambda)\rightarrow B_{n+1}(\lambda)$ \\ 
$b+qL_{n}(\lambda)\mapsto b+qL_{n+1}(\lambda)$%
\end{tabular}
\right. .
\end{equation*}
By identifying $B_{n}(\lambda)$ with its image under $\iota_{n},$ one can
write $B_{\infty}(\lambda)=\cup_{n\geq m}B_{n}(\lambda)$.

\begin{example}
The crystal graph of the vector representation of $U_{q}(g_{B_{\infty}})$ is
identified with the infinite crystal 
\begin{equation*}
B_{B_{\infty}}(1):\cdot\cdot\cdot\rightarrow\overline{n}\overset{n-1}{%
\rightarrow}\overline{n-1}\overset{n-2}{\rightarrow}\cdot\cdot\cdot
\cdot\rightarrow\overline{2}\overset{1}{\rightarrow}\overline{1}\overset{0}{%
\rightarrow}0\overset{0}{\rightarrow}1\overset{1}{\rightarrow}2\cdot
\cdot\cdot\cdot\overset{n-2}{\rightarrow}n-1\overset{n-1}{\rightarrow }%
n\rightarrow\cdot\cdot\cdot
\end{equation*}
\end{example}

\noindent The following proposition is an immediate consequence of the
labelling of the crystals $B_{n}(\lambda)$ obtained in \cite{KN}.

\begin{proposition}
\label{prop_labelcrys}Consider $\lambda\in\mathcal{P}_{m}$. The vertices of $%
B_{\infty}(\lambda)$ are indexed by the tableaux of $\mathbf{T}(\lambda).$
\end{proposition}

\noindent Let $w$ be a word on $\mathcal{X}_{\infty}$. We have seen that
there exists an integer $N$ such that for any integer $k\geq N,$ $P_{k}(w)$
does not depend on $k$ and can be computed without contraction relation
(i.e. the number of boxes in $P(w)$ is equal to the length of $w$). Thus it
makes sense to consider the map 
\begin{equation*}
P_{\infty}:\left\{ 
\begin{tabular}{l}
$\mathcal{X}_{\infty}\rightarrow\mathbf{T}$ \\ 
$w\mapsto\lim_{n\rightarrow\infty}P_{n}(w)$%
\end{tabular}
\right. .
\end{equation*}
Set $G_{\infty}=\underset{l\geq0}{\tbigoplus }B_{\infty}(1)^{\otimes l}$.\
We write $w_{1}\sim w_{2}$ if and only if $w_{1}$ and $w_{2}$ belong to the
same connected component of the crystal $G_{\infty}.\;$

\begin{definition}
The plactic monoids $Pl(A_{\infty}),Pl(B_{\infty}),Pl(C_{\infty})$ and $%
Pl(D_{\infty})$ are the quotient of the free monoids on the alphabets $%
\mathcal{X}_{A_{\infty}},\mathcal{X}_{B_{\infty}},\mathcal{X}_{C_{\infty}},$
and $\mathcal{X}_{D_{\infty}}$ respectively by the relations (\ref{plA}), (%
\ref{plB}), \ref{plC}) and (\ref{plD}).
\end{definition}

\begin{proposition}
\label{prop_P}Consider $w_{1}$ and $w_{2}$ in $\mathcal{X}_{\infty}^{\ast}$%
.\ We have the equivalences 
\begin{equation*}
w_{1}\sim w_{2}\Longleftrightarrow w_{1}\equiv w_{2}\Longleftrightarrow
P_{\infty}(w_{1})=P_{\infty}(w_{2}).
\end{equation*}
\end{proposition}

\begin{proof}
By \cite{Lec} and \cite{lec2} we have for any integer $n\geq1$ the
equivalences 
\begin{equation*}
w_{1}\sim_{n}w_{2}\Longleftrightarrow
w_{1}\equiv_{n}w_{2}\Longleftrightarrow P_{n}(w_{1})=P_{n}(w_{2}).
\end{equation*}
For $n$ sufficiently large we have seen that $P_{n}(w_{1})=P_{\infty}(w_{1})$
and $P_{n}(w_{2})=P_{\infty}(w_{2})$. Moreover there is no contraction
during the computation of $P_{n}(w_{1})$ and $P_{n}(w_{2}).$ Thus for such $%
n $ we have $w_{1}\equiv_{n}w_{2}\Longleftrightarrow w_{1}\equiv w_{2}$.
This yields the equivalence $w_{1}\equiv w_{2}\Longleftrightarrow
P_{\infty}(w_{1})=P_{\infty}(w_{2})$. The equality $P_{\infty}(w_{1})=P_{%
\infty}(w_{2})$ means that for any integer $p\geq n,$ the words $w_{1}$ and $%
w_{2}$ occur at the same place in two isomorphic connected components of $%
B_{p}(1)^{\otimes \ell}$ where $\ell$ is the length of the words $w_{1}$ and 
$w_{2}$. Each connected component $\Gamma$ of $B_{\infty}(1)^{\ell}$ can be
realized as the union of the connected components $\Gamma_{p}$ of $%
B_{p}(1)^{\otimes\ell}$ $p\geq n$ obtained by deleting in $\Gamma$ the
arrows of colors $i\notin I_{p}$ together with the vertices they connect.
Thus we have $P_{\infty}(w_{1})=P_{\infty}(w_{2})\Longleftrightarrow
w_{1}\sim w_{2}$.
\end{proof}

\noindent By the above proposition the plactic classes are labelled by the
tableaux of $\mathbf{T}$.$\;$In the sequel we will write $[w]$ for the
plactic class of the word $w$. We have in particular $\mathrm{w}%
(P_{\infty}(w))\in\lbrack w]$.

\subsection{RS-correspondences and bi-crystal structures\label{subsec_RS}}

By Proposition \ref{prop_labelcrys}, for any positive integer $k,$ the
vertices of $B_{\infty}(k)$ are labelled by the row tableaux of length $k$.\
The combinatorial description of these tableaux depends on the type
considered. More precisely a row tableau $R$ of type $\Delta_{\infty}$ ($%
\Delta=A,B,C,D$) and length $k$ is a row Young diagram of length $k$ filled
from left to right by weakly increasing letters of $\mathcal{X}_{\Delta
_{\infty}}$ and containing at most one letter $0$ when $\Delta_{\infty
}=B_{\infty}$. This implies that there is two kinds of rows of type $%
B_{\infty}$ following they contain the letters $0$ or not.\ Similarly, since
the letters $1$ and $\overline{1}$ are not comparable in $D_{\infty},$ there
is also two kinds of rows of type $D_{\infty}$ following they contain at
least a letter $1$ or not. Accordingly to the convention for the reading of
a tableau, the reading of the row $R$ of length $r$ can be written $\mathrm{w%
}(R)=x_{1}\cdot\cdot\cdot x_{r}$ where $x_{1}\geq\cdot\cdot\cdot\geq x_{r}.$

\noindent Let $k$ be a nonnegative integer and consider $\tau=(\tau
_{1},...,\tau_{k})\in\mathbb{N}^{k}$.\ Set $B_{\infty}(\tau)=B_{\infty}(%
\tau_{1})\otimes\cdot\cdot\cdot\otimes B_{\infty}(\tau_{k})$ and 
\begin{equation*}
\mathcal{W}_{\infty}=\bigoplus_{k=0}^{\infty}\bigoplus_{\tau\in\mathbb{N}%
^{k}}^{{}}B_{\infty}(\tau).
\end{equation*}
Similarly for any integer $n\geq1,$ we define $B_{n}(\tau)=B_{n}(\tau
_{1})\otimes\cdot\cdot\cdot\otimes B_{n}(\tau_{k})$ and 
\begin{equation*}
\mathcal{W}_{n}=\bigoplus_{k=0}^{\infty}\bigoplus_{\tau\in\mathbb{N}%
^{k}}^{{}}B_{n}(\tau).
\end{equation*}
Let $\mathcal{Y}_{\infty}=\{y_{1}<y_{2}<\cdot\cdot\cdot\}$ be a totally
ordered alphabet labelled by the positive integers. A bi-word on $(\mathcal{X%
}_{\infty},\mathcal{Y}_{\infty})$ is a word in the commutative variables $%
\binom{y_{j}}{x_{i}}$ where $x_{i}\in\mathcal{X}_{\infty}$ and $y_{j}\in%
\mathcal{Y}_{\infty}$.\ For types $A_{\infty},B_{\infty}$ and $C_{\infty}$
each biword $W$ can be uniquely written in the form 
\begin{equation}
W=\left( 
\begin{array}{lllllllllllllll}
y_{1} & \cdot & \cdot & y_{1} & y_{2} & \cdot & \cdot & y_{2} & \cdot & \cdot
& \cdot & y_{k} & \cdot & \cdot & y_{k} \\ 
b_{1} & \cdot & \cdot & e_{1} & b_{2} & \cdot & \cdot & e_{2} & \cdot & \cdot
& \cdot & b_{k} & \cdot & \cdot & e_{k}
\end{array}
\right) =\left( 
\begin{array}{l}
w_{y} \\ 
w_{x}
\end{array}
\right)  \label{biwor}
\end{equation}
by ordering first the letters $y_{j}$ next the letters $x_{i}$ so that for
any $p\in\{1,...,k\}$, $b_{p}\cdot\cdot\cdot e_{p}$ is a word on $\mathcal{X}%
_{\infty}$ verifying $b_{p}\geq\cdot\cdot\cdot\geq e_{p}$. For type $%
D_{\infty},$ we can obtain similarly a canonical form for the biwords by
requiring that the letters $\overline{1}$ appear on the left of the letters $%
1$ in each word $b_{p}\cdot\cdot\cdot e_{p}$. Now consider a vertex $b\in%
\mathcal{W}_{\infty}.\;$There exists an integer $k$ and a $k$-uple $%
\tau=(\tau_{1},...,\tau_{k})$ such that $b\in B_{\infty}(\tau)=B_{\infty
}(\tau_{1})\otimes\cdot\cdot\cdot\otimes B_{\infty}(\tau_{k})$.\ Write $%
b=R_{1}\otimes\cdot\cdot\cdot\otimes R_{k}$ with $R_{p}\in
B_{\infty}(\tau_{p})$ for any $p\in\{1,...,k\}.\;$In the sequel we will
identify $b$ with the biword (\ref{biwor}) where for any $p=1,...,k,$ $%
b_{p}\cdot\cdot\cdot e_{p}$ is the reading of the row $R_{p}.$ For types $%
A_{\infty}$ and $C_{\infty},$ each biword can be interpreted as a vertex of $%
\mathcal{W}_{\infty}$ by the above identification.\ This is not true for
types $B_{\infty}$ and $D_{\infty}$ due to the restrictive conditions on
rows.

\noindent Consider $b\in B_{\infty}^{\tau}$.\ By abuse of notation, we write 
$P_{\infty}(b)$ for $P_{\infty}(w_{x})$.\ For any $p\in\{1,...,k\}$ write $%
Q^{(p)}$ for the shape of the tableau $P_{\infty}(b_{1}\cdot\cdot\cdot
e_{p}).\;$The shape $Q^{(p+1)}$ is obtained by adding $\tau_{p+1}$ boxes to $%
Q^{(p)},$ thus we can consider that $Q^{(p)}\subset Q^{(p+1)}$.\ Denote by $%
Q_{\infty}(b)$ the tableau obtained by filling $Q^{(1)}$ with letters $%
y_{1}, $ next each $Q^{(p+1)}\backslash Q^{(p)}$ with letters $y_{p+1}.$

\begin{lemma}
\label{lem_sst}$Q_{\infty}(b)$ is a semistandard tableau on $\mathcal{Y}%
_{\infty}$.
\end{lemma}

\begin{proof}
By induction on $k$, it suffices to prove that the insertion of the row $R$
in a tableau $T$ gives a tableau $T^{\prime}$ such that $T^{\prime}%
\backslash T$ is an horizontal strip skew tableau (that is has no boxes in
the same column). There exists an integer $n$ such that $R,T$ and $%
T^{\prime} $ are tableaux of $\mathbf{T}^{n}$ verifying $T^{\prime}=P_{n}(%
\mathrm{w}(T)\mathrm{w}(R)).$ Let $\lambda$ be the shape of $T$ and $r$ the
length of $R$. Write $B(\mathrm{w}(T)\otimes\mathrm{w}(R))$ for the
connected component of $B_{n}(\lambda )\otimes B(1^{r})$ containing $\mathrm{%
w}(T)\otimes\mathrm{w}(R).\;$The shape $\nu$ of $T^{\prime}$ is given by the
weight of the highest weight vertex in $B(\mathrm{w}(T)\otimes\mathrm{w}(R)).
$ By using Lemma \ref{lem_plu_hp} and the fact that $\left| \nu\right|
=\left| \lambda\right| +r,$ one easily obtains that $\nu/\lambda$ is an
horizontal strip.
\end{proof}

\noindent In the sequel we will denote by $\mathrm{SST}$ the set of
semistandard tableaux on the alphabet $\mathcal{Y}_{\infty}$.

\bigskip

\noindent\textbf{Remarks:}

\noindent$\mathrm{(i)}\mathbf{:}$ Given a word $w,$ its $Q$-symbol depends
on the type considered.\ For example if we take $b=1\otimes\overline{1}%
\otimes1\otimes\overline{1}\in B^{C_{\infty}}(1,1,1,1)\cap
B^{D_{\infty}}(1,1,1,1)$, we will have 
\begin{equation*}
\left( P_{C_{\infty}}(b),Q_{C_{\infty}}(b)\right) =\left( 
\begin{tabular}{|l|l|}
\hline
$\mathtt{\bar{2}}$ & $\mathtt{1}$ \\ \hline
$\mathtt{\bar{1}}$ & $\mathtt{2}$ \\ \hline
\end{tabular}
, 
\begin{tabular}{|l|l|}
\hline
$\mathtt{1}$ & $\mathtt{2}$ \\ \hline
$\mathtt{3}$ & $\mathtt{4}$ \\ \hline
\end{tabular}
\right) \text{ and }\left( P_{D_{\infty}}(b),Q_{D_{\infty}}(b)\right)
=\left( 
\begin{tabular}{|l|}
\hline
$\mathtt{1}$ \\ \hline
$\mathtt{\bar{1}}$ \\ \hline
$\mathtt{1}$ \\ \hline
$\mathtt{\bar{1}}$ \\ \hline
\end{tabular}
\text{ } 
\begin{tabular}{|l|}
\hline
$\mathtt{1}$ \\ \hline
$\mathtt{2}$ \\ \hline
$\mathtt{3}$ \\ \hline
$\mathtt{4}$ \\ \hline
\end{tabular}
\right) .
\end{equation*}

\noindent$\mathrm{(ii)}\mathbf{:}$ $Q_{\infty}(b)$ can also be considered as
the sequence of the shapes of the tableaux $P_{\infty}(b_{1}\cdot\cdot\cdot
e_{p}),$ $p\in\{1,...,k\}.$

\noindent$\mathrm{(iii)}\mathbf{:}$ It follows immediately from their
definition that $P_{\infty}(b)$ and $Q_{\infty}(b)$ have the same shape.

\begin{proposition}
\label{prop_Q}Consider a semistandard tableau $t$\ of shape $\lambda$ and
set $B_{\infty}(t)=\{b\in\mathcal{W}_{\infty}\mid Q_{\infty}(b)=t\}.\;$Then $%
B_{\infty}(t)$ is a connected component of $\mathcal{W}_{\infty}$ isomorphic
to $B_{\infty}(\lambda)$.
\end{proposition}

\begin{proof}
By the previous lemma, we know that $Q_{\infty}(b)$ is a semi-standard
tableau.\ It follows from \cite{Lec} and \cite{lec2} that for any integer $%
n\geq1,$ $B_{n}(t)=B_{\infty}(t)\cap\mathcal{W}_{n}$ is a connected
component isomorphic to $B_{n}(\lambda).$ Since $B_{\infty}(t)=\cup_{n%
\geq1}B_{n}(t)$ and $B_{\infty}(\lambda)=\cup_{n\geq1}B_{n}(\lambda)$ this
gives the proposition.
\end{proof}

\bigskip

\noindent For each $\Delta_{\infty}\in\{A_{\infty},B_{\infty},C_{\infty
},D_{\infty}\},$ set $\mathcal{E}_{\Delta_{\infty}}=\{(P,Q)\mid P$ and $Q$
are tableaux with the same shape, respectively of type $\Delta_{\infty}$ and
semistandard\ on $\mathcal{Y}_{\infty}\}.$ The following theorem can be
regarded as a generalization of the Robinson-Schensted correspondence on
biwords:

\begin{theorem}
\label{th_RS}The map 
\begin{equation*}
\mathrm{RS}:\left\{ 
\begin{tabular}{l}
$\mathcal{W}_{\infty}\rightarrow\mathcal{E}_{\infty}$ \\ 
$b\mapsto(P_{\infty}(b),Q_{\infty}(b))$%
\end{tabular}
\right.
\end{equation*}
is a one to one correspondence.
\end{theorem}

\begin{proof}
Consider $b_{1}$ and $b_{2}$ in $\mathcal{W}_{\infty}$ such that $\mathrm{RS}%
(b_{1})=\mathrm{RS}(b_{2})$. Then by Propositions \ref{prop_P} and \ref
{prop_Q}, $b_{1}$ and $b_{2}$ occur at the same place in the same connected
component of $\mathcal{W}_{\infty}$. Hence $b_{1}=b_{2}$ and $\mathrm{RS}$
is an injective map. Now consider $(P,Q)\in\mathcal{E}_{\infty}$. There
exist an integer $n\geq1$ such that $P\in\mathbf{T}_{n}.$ Then by \cite{Lec}
and \cite{lec2}, there exists $b\in\mathcal{W}_{n}\subset \mathcal{W}%
_{\infty}$ such that $P_{n}(b)=P$ and $Q_{n}(b)=Q$. In this case we have $%
P_{\infty}(b)=P$ and $Q_{\infty}(b)=Q(b)$, thus $\mathrm{RS}$ is an
surjective map.
\end{proof}

\bigskip

\noindent\textbf{Remarks: }

\noindent$\mathrm{(i)}\mathbf{:}$ Let $b$ be a vertex of $B_{n}^{\tau}$ and
denote by $Q_{n}(b)=(Q_{n}^{(1)},...,Q_{n}^{(k)})$ the sequence of shapes of
the tableaux $P_{n}(b_{1}\cdot\cdot\cdot e_{p}),$ $p\in\{1,...,k\}.$ In
general $Q_{n}^{(p+1)}\backslash Q_{n}^{(p)}$ is not an horizontal strip for
the number of boxes in $Q_{n}^{(p+1)}$ can be strictly less than that of $%
Q_{n}^{(p)}$. This implies in particular that the map $b\longmapsto\left(
P_{C_{n}}(b),Q_{C_{n}}(b)\right) $ defined on $\mathcal{W}_{n}$ is not a
one-to-one correspondence.\ 

\noindent$\mathrm{(ii)}\mathbf{:}$ When any two consecutive shapes in $%
Q_{n}(b)$ differ by at most one box, one says that $Q_{n}(b)$ is an
oscillating tableau.\ By (\cite{Lec} and \cite{lec2}), the map 
\begin{equation*}
\mathrm{rs}_{n}:\left\{ 
\begin{array}{l}
\bigoplus_{k=0}^{\infty}B_{n}(1^{k})\rightarrow\mathcal{O}_{n} \\ 
b\mapsto(P_{n}(b),Q_{n}(b))
\end{array}
\right.
\end{equation*}
where $\mathcal{O}_{\Delta_{n}}=\{(P,Q)\mid P$ is a tableau and $Q$ an
oscillating tableau\ of type $\Delta_{n}\}$ is then a one-to-one
correspondence. This difference with the situation described in $\mathrm{(i)}
$ follows from the fact that the decomposition $V_{n}(\lambda)\otimes
V_{n}(1^{k})$ is not free multiplicity for types $B_{n},C_{n}$ and $D_{n}$
when $k>1$.

\noindent$\mathrm{(iii)}\mathbf{:}$ The map 
\begin{equation*}
\mathrm{rs}_{\infty}:\left\{ 
\begin{array}{l}
\bigoplus_{k=0}^{\infty}B_{\infty}(1^{k})\rightarrow\mathcal{O}_{\infty} \\ 
b\mapsto(P_{\infty}(b),Q_{\infty}(b))
\end{array}
\right.
\end{equation*}
where $\mathcal{O}_{\Delta_{\infty}}=\{(P,Q)\mid P$ and $Q$ are tableaux \
respectively of type $\Delta_{\infty}$ and standard on $\mathcal{Y}_{\infty
} $ with the same shape$\}$ is also a one-to-one correspondence. Thus for
any word $w,$ the plactic class $[w]$ contains exactly $n_{\lambda}$ words
(see Theorem \ref{th_decovec}) where $\lambda$ is the shape of $P(w).$ In
particular, it is finite.

\bigskip

\noindent By Proposition \ref{prop_labelcrys}, the set $\mathbf{T}(\lambda)$
of tableaux of type $\Delta_{\infty}=B_{\infty},C_{\infty},D_{\infty}$ can
be endowed with the structure of a $U_{q}(\frak{g}_{\Delta_{\infty}})$%
-crystal so that for any two tableaux $T$ and $T^{\prime}$ in $\mathbf{T}%
(\lambda)$ and any Kashiwara operator $\widetilde{k}_{i}=\widetilde{e}_{i},%
\widetilde{f}_{i},$ $i\geq0$%
\begin{equation*}
\widetilde{k}_{i}(T)=T^{\prime}\Longleftrightarrow\widetilde{k}_{i}(\mathrm{w%
}(T))=\mathrm{w}(T^{\prime})\text{ in }B_{\Delta_{\infty}}(1)^{\otimes\left|
\lambda\right| }.
\end{equation*}
We proceed similarly to define a structure of $U_{q}(\frak{g}_{A_{\infty}})$%
-crystal on the set $\mathbf{SST}(\lambda)$ of semistandard tableaux of
shape $\lambda$ on $\mathcal{Y}_{\infty}$. For any two tableaux $t$ and $%
t^{\prime}$ in $\mathbf{SST}(\lambda)$ and any Kashiwara operator $%
\widetilde{K}_{j}=\widetilde{E}_{j},\widetilde{F}_{j}$, $j\geq1$ we set 
\begin{equation*}
\widetilde{K}_{j}(t)=t^{\prime}\Longleftrightarrow\widetilde{K}_{j}(\mathrm{w%
}(t))=\mathrm{w}(t^{\prime})\text{ in }B_{A_{\infty}}(1)^{\otimes \left|
\lambda\right| }.
\end{equation*}
By definition $\mathcal{W}_{\infty}$ is a $U_{q}(\frak{g}_{\infty})$%
-crystal. We can also endow it with the structure of a $U_{q}(\frak{g}%
_{A_{\infty}})$-crystal by setting 
\begin{equation*}
\widetilde{K}_{j}(b)=\mathrm{RS}^{-1}(P_{\infty}(b),\widetilde{K}%
_{j}(Q_{\infty}(b)))
\end{equation*}

\begin{theorem}
\label{th_bicrysl}$\mathcal{W}_{\infty}$ has the structure of a bi-crystal
for the pair $(U_{q}(\frak{g}_{\infty}),U_{q}(\frak{g}_{A_{\infty}}))$, that
is for any Kashiwara operators $\widetilde{k}_{i},i\in I$ and $\widetilde{K}%
_{j},j\geq1$ the following diagram commutes: 
\begin{equation*}
\begin{tabular}{lll}
$\ \ \mathcal{W}_{\infty}$ & $\overset{\widetilde{K}_{j}}{\rightarrow}$ & $%
\mathcal{W}_{\infty}$ \\ 
$\widetilde{k}_{i}\downarrow$ &  & $\ \downarrow$ $\ \widetilde{k}_{i}$ \\ 
$\ \ \mathcal{W}_{\infty}$ & $\overset{\widetilde{K}_{j}}{\rightarrow}$ & $%
\mathcal{W}_{\infty}$%
\end{tabular}
.
\end{equation*}
\end{theorem}

\begin{proof}
Consider $b\in\mathcal{W}_{\infty}$ and set $\mathrm{RS}(b)=(P,Q)$.\ Then $%
\mathrm{RS}(\widetilde{k}_{i}(b))=(\widetilde{k}_{i}(P),Q)$ for $b$ and $%
\widetilde{k}_{i}(b)$ are in the same $U_{q}(\frak{g}_{\infty})$-connected
component (see Proposition \ref{prop_Q}). This gives $\mathrm{RS}(\widetilde{%
K}_{j}(\widetilde{k}_{i}(b)))=(\widetilde{k}_{i}(P),\widetilde {K}_{j}(Q))$.
One the other hand we have $\mathrm{RS}(\widetilde{K}_{j}(b))=(P,\widetilde{K%
}_{j}(Q))$. Moreover $\mathrm{RS}(\widetilde{k}_{i}(\widetilde{K}_{j}(b)))=(%
\widetilde{k}_{i}(P),\widetilde{K}_{j}(Q))$ since $\widetilde{K}_{j}(b)$ and 
$\widetilde{k}_{i}(\widetilde{K}_{j}(b))$ are in the same $U_{q}(\frak{g}%
_{\infty})$-connected component. We obtain $\widetilde{K}_{j}(\widetilde{k}%
_{i}(b))=\widetilde{k}_{i}(\widetilde{K}_{j}(b))$ because the map $\mathrm{RS%
}$ is a one-to-one correspondence
\end{proof}

\bigskip

\noindent The RS-type correspondences and the bi-crystal structures, we have
obtained in this paragraph can be referred as the ``symmetric case'' by
analogy to the situation in type $A$ where, for each $k\in\mathbb{N}$, the
crystal $B_{\infty}(k)$ is that of the $k$-th symmetric power of the vector
representation. It is also possible to define a RS-type correspondence and a
bi-crystal structure in the ``anti-symmetric case'', that is when $\widehat{%
\mathcal{W}}_{\infty}$ is replaced by the crystal 
\begin{equation*}
\widehat{\mathcal{W}}_{\infty}=\bigoplus_{k=0}^{\infty}\bigoplus_{\tau \in%
\mathbb{N}^{k}}^{{}}\widehat{B}_{\infty}(\tau)
\end{equation*}
where for any $\tau=(\tau_{1},...,\tau_{k})\in\mathbb{N}^{k},$ $\widehat
{B}_{\infty}(\tau)=B_{\infty}(1^{\tau_{1}})\otimes\cdot\cdot\cdot\otimes
B_{\infty}(1^{\tau_{k}}).$ By Proposition \ref{prop_labelcrys}, the vertices
of $B_{\infty}(1^{k}),k\in\mathbb{N}$ are labelled by the columns of height $%
k$ (see (\ref{col})). Consider $b=C_{1}\otimes\cdot\cdot\cdot\otimes C_{k}\in%
\widehat{B}_{\infty}(\tau)$. One associates to $b$ the tableau $\widehat
{P}_{\infty}(b)=P_{\infty}(\mathrm{w}(C_{1})\cdot\cdot\cdot\mathrm{w}%
(C_{k})).$ Write $(Q_{n}^{(1)},...,Q_{n}^{(k)})$ for the sequence of shapes
of the tableaux $P_{\infty}(\mathrm{w}(C_{1})\cdot\cdot\cdot\mathrm{w}%
(C_{p})),$ $p\in\{1,...,k\}.$ One proves that $Q_{n}^{(p+1)}\backslash
Q_{n}^{(p)}$ is a \textit{vertical} strip (that is, has no boxes in the same
row). Let $\lambda$ be the shape of $\widehat{P}_{\infty}(b).$ Denote by $%
\widehat{Q}_{\infty}(b)$ the tableau obtained by filling $Q^{(1)}$ with
letters $y_{1},$ next each $Q^{(p+1)}\backslash Q^{(p)}$ with letters $%
y_{p+1}$. The rows of $\widehat {Q}_{\infty}(b)$ (read from left to right)
strictly increase and its columns (read from top to bottom) weakly increase.
This means that $\widehat {Q}_{\infty}(b)$ can be regarded as a semistandard
tableau of shape $\lambda^{\prime},$ the conjugate partition of $\lambda.$

\noindent For each $\Delta_{\infty}\in\{A_{\infty},B_{\infty},C_{\infty
},D_{\infty}\},$ set $\widehat{\mathcal{E}}_{\Delta_{\infty}}=\{(P,Q)\mid P$
and $Q$ are tableaux with conjugate shapes, respectively of type $%
\Delta_{\infty}$ and semistandard\ on $\mathcal{Y}_{\infty}\}.$ By
proceeding as for Theorem \ref{th_RS}, we obtain:

\begin{theorem}
\label{th_bicrys2}The map 
\begin{equation*}
\widehat{\mathrm{RS}}:\left\{ 
\begin{tabular}{l}
$\widehat{\mathcal{W}}_{\infty}\rightarrow\widehat{\mathcal{E}}_{\infty}$ \\ 
$b\mapsto(\widehat{P}_{\infty}(b),\widehat{Q}_{\infty}(b))$%
\end{tabular}
\right.
\end{equation*}
is a one to one correspondence.
\end{theorem}

\noindent This also permits to endow $\widehat{\mathcal{W}}_{\infty}$ with a
structure of a $U_{q}(\frak{g}_{A_{\infty}})$-crystal by setting 
\begin{equation*}
\widetilde{K}_{j}(b)=\widehat{\mathrm{RS}}^{-1}(P_{\infty}(b),\widetilde
{K}_{j}(Q_{\infty}(b)))
\end{equation*}
for any Kashiwara operator $\widetilde{K}_{j},j\geq1$ associated to $U_{q}(%
\frak{g}_{A_{\infty}}).$ Similarly to Theorem \ref{th_bicrysl}, we derive

\begin{theorem}
$\widehat{\mathcal{W}}_{\infty}$ has the structure of a bi-crystal for the
pair $(U_{q}(\frak{g}_{\infty}),U_{q}(\frak{g}_{A_{\infty}})).$
\end{theorem}

\section{A plactic algebra for infinite root systems}

\subsection{Plactic Schur functions}

\noindent Denote by $\mathbb{Z}[\mathcal{W}_{\infty}]$ the free $\mathbb{Z}$%
-algebra of basis $\mathcal{W}_{\infty}$ where the multiplication is defined
by concatenation of tensors. More precisely, for any $\tau=(\tau_{1},...,%
\tau_{k})\in\mathbb{N}^{k},\tau=(\tau_{1}^{\prime},...,\tau_{k}^{\prime })\in%
\mathbb{N}^{k^{\prime}}$ and $b=b_{1}\otimes\cdot\cdot\cdot b_{r}\in
B_{\infty}(\tau),b^{\prime}=b_{1}^{\prime}\otimes\cdot\cdot\cdot
b_{r}^{\prime}\in B_{\infty}(\tau^{\prime})$ we have 
\begin{equation*}
b\cdot b^{\prime}=b_{1}\otimes\cdot\cdot\cdot\otimes b_{r}\otimes
b_{1}^{\prime}\otimes\cdot\cdot\cdot\otimes b_{r}^{\prime}.
\end{equation*}
Consider a partition $\lambda$ and $t$ a semistandard tableau of shape $%
\lambda$ on the alphabet $\mathcal{Y}_{\infty}$.\ Denote by $B_{\infty}(t)$
the connected component of $\mathcal{W}_{\infty}$ associated to $t$ by
Proposition \ref{prop_Q}. We define the free Schur function associated to $t$
on $\mathbb{Z}[\mathcal{W}_{\infty}]$ by 
\begin{equation*}
\mathbf{S}_{t}=\sum_{Q_{\infty}(b)=t}b=\sum_{b\in B_{\infty}(t)}b.
\end{equation*}
For any $b=b_{1}\otimes\cdot\cdot\cdot\otimes b_{r}\in B_{\infty}(\tau),$
set $\mathrm{w}(b)=\mathrm{w}(b_{1})\cdot\cdot\cdot\mathrm{w}(b_{r}).\;$This
is equivalent to say that $\mathrm{w}(b)=w_{x}$ when $b$ is considered as
the biword $\binom{w_{y}}{w_{x}}.$ The plactic algebra $\mathbb{Z}%
[Pl(\infty)]$ is the quotient of $\mathbb{Z}[\mathcal{W}_{\infty}]$ by the
relations 
\begin{equation*}
b=b^{\prime}\Longleftrightarrow\mathrm{w}(b)\equiv\mathrm{w}(b^{\prime}).
\end{equation*}
Denote by $\frak{P}_{\infty}$ the canonical projection from $\mathbb{Z}[%
\mathcal{W}_{\infty}]$ to the plactic algebra $\mathbb{Z}[Pl(\infty)].$ For
any $b\in\mathcal{W}_{\infty}$, we thus have $\frak{P}_{\infty}(b)=P(\mathrm{%
w}(b)).$ By Proposition \ref{prop_P}, two vertices $b_{1}$ and $b_{2}$ in $%
\mathcal{W}_{\infty}$ occur at the same place in two isomorphic components
if and only if $\frak{P}_{\infty}(b_{1})=\frak{P}_{\infty}(b_{2})$. Since
the plactic classes are labelled by the tableaux, $\{T\mid T\in\mathbf{T}\}$
is a basis of $\mathbb{Z}[Pl(\infty)].$ By Propositions \ref{prop_P} and \ref
{prop_Q} 
\begin{equation*}
\frak{P}_{\infty}(\mathbf{S}_{t})=\sum_{T\in B_{\infty}(\lambda)}T
\end{equation*}
depends only of the shape $\lambda.\;$The right member of the above equality
is called the plactic Schur function and is denoted by $\mathbf{S}_{\lambda
}^{\infty}$.

\begin{theorem}
For each type $A_{\infty},B_{\infty},C_{\infty},D_{\infty},$ the plactic
Schur functions span a commutative subalgebra of $\mathbb{Z}[Pl(\infty)]$.\
Moreover given two partitions $\lambda$ and $\mu$ we have 
\begin{equation*}
\mathbf{S}_{\lambda}^{\infty}\cdot\mathbf{S}_{\mu}^{\infty}=\sum_{\nu \in%
\mathcal{P}}c_{\lambda,\mu}^{\nu}\mathbf{S}_{\nu}^{\infty}
\end{equation*}
where $c_{\lambda,\mu}^{\nu}$ is a Littlewood-Richardson coefficient, thus
does not depend on the type considered.
\end{theorem}

\begin{proof}
In the plactic algebra $\mathbb{Z}[Pl(\infty)],$ we have 
\begin{equation*}
\mathbf{S}_{\lambda}^{\infty}\cdot\mathbf{S}_{\mu}^{\infty}=\sum_{T^{\prime
}\in B_{\infty}(\lambda)}\sum_{T^{^{\prime\prime}}\in B_{\infty}(\mu
)}T^{\prime}\cdot T^{^{\prime\prime}}=\sum_{T^{\prime}\otimes
T^{\prime\prime }\in B_{\infty}(\lambda)\otimes B_{\infty}(\mu)}T\cdot
T^{\prime\prime}=\sum_{\nu}\sum_{T\in B_{\infty}(\nu)}m(T)T
\end{equation*}
where $m(T)$ is the number of vertices $T^{\prime}\otimes
T^{\prime\prime}\in B_{\infty}(\lambda)\otimes B_{\infty}(\mu)$ such that $%
T^{\prime}\cdot T^{^{\prime\prime}}=T$ in $\mathbb{Z}[Pl(\infty)]$. Let us
write $\emptyset$ for the empty crystal.\ By assertion $2$ of Corollary \ref
{dec_tens}, Proposition \ref{prop_sdpt} and Theorem \ref{th_iso_bc}, there
is a crystal isomorphism 
\begin{equation*}
\Upsilon_{\lambda,\mu}:B_{\infty}(\lambda)\otimes B_{\infty}(\mu )\overset{%
\simeq}{\rightarrow}\bigoplus_{\nu}\bigoplus_{r=1}^{c_{\lambda,\mu
}^{\nu}}B_{\infty}(\nu)\otimes\emptyset^{\otimes r}.
\end{equation*}
Since the plactic relations express crystal isomorphisms, for $T,T^{\prime
},T^{\prime\prime}$ respectively in $B_{\infty}(\nu),B_{\infty}(\lambda)$
and $B_{\infty}(\mu)$, we have the equivalences 
\begin{equation*}
T^{\prime}\cdot T^{^{\prime\prime}}=T\Longleftrightarrow\mathrm{w}(T^{\prime
})\mathrm{w}(T^{\prime\prime})\equiv\mathrm{w}(T)\Longleftrightarrow\exists
r\in\{1,...,c_{\lambda,\mu}^{\nu}\}\mid\Upsilon_{\lambda,\mu}(T^{\prime
}\otimes T^{\prime\prime})=T\otimes\emptyset^{\otimes r}.
\end{equation*}
This means that $m(T)=c_{\lambda,\mu}^{n}$ for any $T\in B_{\infty}(\nu)$.
Hence we obtain 
\begin{equation*}
\mathbf{S}_{\lambda}^{\infty}\cdot\mathbf{S}_{\mu}^{\infty}=\sum_{\nu
}c_{\lambda,\mu}^{\nu}\sum_{T\in B_{\infty}(\nu)}T=\sum_{\nu}c_{\lambda,\mu
}^{\nu}\mathbf{S}_{\nu}^{\infty}.
\end{equation*}
\end{proof}

\noindent We deduce immediately the following corollary:

\begin{corollary}
The four algebras generated by the plactic Schur functions of types $%
A_{\infty},B_{\infty},C_{\infty},D_{\infty}$ are distinct subalgebras of $%
\mathbb{Z}[Pl(\infty)]$ all isomorphic to the algebra of symmetric functions.
\end{corollary}

\subsection{Cauchy-types identities}

Let $\mathcal{A}[[X,Y]]$ be the $\mathbb{Z}$ algebra of formal series in the
noncommutative variables $x_{i}\in\mathcal{X}_{\infty}$ and $y_{j}\in%
\mathcal{Y}_{\infty}$ submit to the relations 
\begin{equation}
\left\{ 
\begin{tabular}{l}
$x_{\overline{i}}x_{i}=x_{i+1}x_{\overline{i+1}}\text{ for any }i\geq1\text{%
, }x_{i^{\prime}}x_{i}=x_{i}x_{i^{\prime}}\text{ for any }(i,i^{\prime})%
\text{ with }i^{\prime}\neq-i$ \\ 
$y_{j}y_{j^{\prime}}=y_{j^{\prime}}y_{j}$ for any $(j,j^{\prime})$ and $%
x_{i}y_{j}=y_{j}x_{i}$ for any $(i,j)$%
\end{tabular}
\right.  \label{rel}
\end{equation}
By using Relations (\ref{rel}), each monomial in the variables $x_{i}$ can
be uniquely written in the form 
\begin{equation*}
x^{\alpha}=x_{n}^{\alpha_{n}}\cdot\cdot\cdot
x_{1}^{\alpha_{1}}x_{0}^{\alpha_{0}}x_{\overline{1}}^{\alpha_{\overline{1}%
}}\cdot\cdot\cdot x_{\overline{m}}^{\alpha_{\overline{m}}}
\end{equation*}
where the $\alpha_{i}$'s are nonnegative integers such that $%
\alpha_{n}\alpha_{\overline{m}}\neq0$.

\noindent To each word $w=x_{1}\cdot\cdot\cdot x_{\ell}$ on $\mathcal{X}%
_{\infty},$ we associated the monomial $x^{w}=x_{1}\cdot\cdot\cdot
x_{\ell}\in\mathcal{A}[[X]]$. By (\ref{rel}) 
\begin{equation}
w_{1}\equiv w_{2}\Longrightarrow x^{w_{1}}=x^{w_{2}}  \label{compat}
\end{equation}
for $w_{1}$ and $w_{2}$ two words on $\mathcal{X}_{\infty}$.

\bigskip

\noindent Given $x\in\mathcal{X}_{\infty}$ and $y\in\mathcal{Y}_{\infty}$ we
have $\frac{1}{1-xy}\in\mathcal{A}[[X,Y]]$. Moreover in $\mathcal{A}[[X,Y]]$
we can write 
\begin{align*}
\frac{1}{1-x_{\overline{i}}y}\frac{1}{1-x_{i}y^{\prime}} & =\frac
{1}{1-x_{i+1}y}\frac{1}{1-x_{\overline{i+1}}y^{\prime}}\text{ for any }i\geq1
\\
\frac{1}{1-x_{i}y}\frac{1}{1-x_{j}y^{\prime}} & =\frac{1}{1-x_{j}y^{\prime}}%
\frac{1}{1-x_{i}y}\text{ if }j\neq-i\text{ or }i=j=0.
\end{align*}
For any partition $\lambda,$ denote by $s_{\lambda}(Y)$ the infinite Schur
function of type $A$ in the variables $y\in\mathcal{Y}_{\infty}$. For $%
\Delta=A,B,C,D,$ set 
\begin{equation*}
S_{\lambda}^{\Delta_{\infty}}(X)=\sum_{T\in\mathbf{T}(\lambda)}x^{\mathrm{w}%
(T)}\in\mathcal{A}[[X]].
\end{equation*}
Note that for types $B_{\infty},C_{\infty}$ and $D_{\infty},$ it is
impossible to specialize $x_{\overline{i}}=\frac{1}{x_{i}}$ in $%
S_{\lambda}^{\Delta_{\infty}}(X)$ because $\mathbf{T}(\lambda)$ contains an
infinite number of tableaux with the same weight. Thus $S_{\lambda}^{%
\Delta_{\infty}}(X)$ can not be interpreted as an ordinary character for
classical Lie algebras.

\noindent For each type $\Delta_{\infty}$ ($\Delta=A,B,C,D$) set 
\begin{equation*}
K_{\Delta_{\infty}}(X,Y)=\left\{ 
\begin{array}{l}
\prod_{i=1}^{\infty}\prod_{j=1}^{\infty}\frac{1}{1-x_{i}y_{j}}\text{ if }%
\Delta=A\vspace*{0.15cm} \\ 
\prod_{i=1}^{\infty}\prod_{j=1}^{\infty}\frac{1+x_{0}}{(1-x_{i}y_{j})(1-x_{%
\overline{i}}y_{j})}\text{ if }\Delta=B\vspace*{0.15cm} \\ 
\prod_{i=1}^{\infty}\prod_{j=1}^{\infty}\frac{1}{(1-x_{i}y_{j})(1-x_{%
\overline {i}}y_{j})}\text{ if }\Delta=C\vspace*{0.15cm} \\ 
\prod_{i=2}^{\infty}\prod_{j=1}^{\infty}\frac{1}{(1-x_{i}y_{j})(1-x_{%
\overline {i}}y_{j})}\left( \frac{1}{1-x_{1}y_{j}}+\frac{1}{1-x_{\overline{1}%
}y_{j}}\right) \text{ if }\Delta=D
\end{array}
\right.
\end{equation*}

\begin{theorem}
\label{the_cauchy}For $\Delta=A,B,C,D,$ we have the following Cauchy-type
identities in $\mathcal{A}[[X,Y]]$ 
\begin{equation}
K_{\Delta_{\infty}}(X,Y)=\sum_{\lambda\in\mathcal{P}}S_{\lambda}^{\Delta_{%
\infty}}(X)s_{\lambda}(Y)  \label{cauchy}
\end{equation}
\end{theorem}

\begin{proof}
To each biword $W=\binom{w_{y}}{w_{x}}\in\mathcal{W}_{\infty}$ (see (\ref
{biwor})) such that $w_{x}=x_{i_{1}}\cdot\cdot\cdot x_{i_{r}}$ and $%
w_{y}=y_{i_{1}}\cdot\cdot\cdot y_{i_{r}}$ we associate the monomial $%
x^{w_{x}}y^{w_{y}}=x_{i_{1}}\cdot\cdot\cdot
x_{i_{r}}y_{i_{1}}\cdot\cdot\cdot y_{i_{r}}\in$ $\mathcal{A}[[X,Y]].$ From
the description of the row tableaux given in \ref{subsec_RS} and by using
Relations (\ref{rel}), we obtain in $\mathcal{A}[[X,Y]]$ 
\begin{equation*}
K_{\Delta_{\infty}}(X,Y)=\sum_{W\in\mathcal{W}_{\infty}}x^{w_{x}}y^{w_{y}}.
\end{equation*}
By Theorem \ref{th_RS}, the map $\mathrm{RS}$ is a one-to-one correspondence
between the biwords $W\in\mathcal{W}_{\infty}$ and the pair $(T,t)$ where $T$
and $t$ are respectively a tableau on $\mathcal{X}_{\infty}$ and a
semistandard tableau on $\mathcal{Y}$ with the same shape. Moreover we have $%
x^{w_{x}}=x^{\mathrm{w}(T)}$ and $y^{w_{y}}=y^{\mathrm{w}(t)}$ (see (\ref
{compat})) in $\mathcal{A}[[X,Y]]$.\ This gives 
\begin{equation*}
K_{\Delta_{\infty}}(X,Y)=\sum_{\lambda}\sum_{(T,t)\in\mathbf{T}(\lambda
)\times\mathbf{SST}(\lambda)}x^{\mathrm{w}(T)}y^{\mathrm{w}%
(t)}=\sum_{\lambda }\sum_{T\in\mathbf{T}(\lambda)}x^{\mathrm{w}(T)}\sum_{t\in%
\mathbf{SST}(\lambda)}y^{\mathrm{w}(t)}=\sum_{\lambda\in\mathcal{P}%
}S_{\lambda}^{\Delta_{\infty}}(X)s_{\lambda}(Y).
\end{equation*}
\end{proof}

\newpage 

\noindent\textbf{Remarks:}

\noindent$\mathrm{(i)}\mathbf{:}$ It is possible to derive Cauchy identities
for a finite alphabet $\mathcal{Y}_{n}=\{y_{1}<\cdot\cdot\cdot<y_{n}\}$ by
specializing the variables $y_{i},i>n$ to $0$ in (\ref{cauchy}).

\noindent$\mathrm{(ii)}\mathbf{:}$ Such a specialization is not possible for
the variables $x_{i}$ due to the relations $x_{\overline{i}}x_{i}=x_{i+1}x_{%
\overline{i+1}}$.\ This means that, except for type $A,$ Cauchy types
identities (\ref{cauchy}) make only sense in infinite rank case. This is
related to the fact that the Robinson-Schensted correspondence of Theorem 
\ref{th_RS} cannot be restricted to types $B_{n},C_{n}$ and $D_{n}$ (see
Remark $\mathrm{(iii)}$ after Theorem \ref{th_RS}).

\noindent$\mathrm{(iii)}\mathbf{:}$ The previous arguments shows that it is
not possible to recover the Littelwood formulas \cite{littlw} from (\ref
{cauchy}). We have find no combinatorial interpretation of Littlewood's
identities in terms of a RS-type correspondence.

\section{Combinatorial question}

\noindent The bi-crystal structure we have obtained in Theorems \ref
{th_bicrysl} and \ref{th_bicrys2} are defined implicitly by using the
RS-type correspondences on $\mathcal{W}_{\infty}$ and $\widehat{\mathcal{W}}%
_{\infty}$. By using the results obtained in \cite{DC}, one can easily
deduce that the $(U_{q}(\frak{g}_{A_{\infty}}),U_{q}(\frak{g}_{A_{\infty}}))$
bi-crystals (that is, corresponding to the case $\Delta=A$) can be naturally
labelled by certain infinite arrays making the symmetry between the two
crystal structures apparent. It would be interesting to define analogues of
these arrays for the $(U_{q}(\frak{g}_{\infty}),U_{q}(\frak{g}%
_{A_{\infty}})) $ bi-crystals of Theorems \ref{th_bicrysl} and \ref
{th_bicrys2} when $\frak{g}_{\infty}\neq\frak{g}_{A_{\infty}}$.

\end{document}